Integral Representations for Interface Problems for the Barenblatt-Sobolev-Galpern Pseudoparabolic Equation

Andreas Chatziafratis [1,*] , Sergey A. Rukolaine [2] , Elias C. Aifantis [3-6]

[1] Division of Applied Analysis & PDE, Department of Mathematics, National & Kapodistrian University of Athens
[2] Ioffe Institute, 26 Polytekhnicheskaya, St. Petersburg 194021, Russia
[3] Laboratory of Mechanics and Materials, School of Civil Engineering, Aristotle University of Thessaloniki
[4] Institute of Plasma Physics & Lasers, School of Electrical & Electronic Engineering, Hellenic Mediterranean University
[5] College of Engineering, Michigan Technological University, USA
[6] Friedrich-Alexander University of Erlangen-Nuremberg, Germany

**Abstract.** We obtain novel integral representations, expressed as contour integrals in the complex Fourier plane, for the solution of fully nonhomogeneous initial-boundary-value as well as interface problems for the Barenblatt-Zheltov-Kochina pseudo-parabolic equation of Sobolev-Galpern type formulated on the real line, half-line and finite interval. This fundamental partial differential equation (PDE) of mathematical physics emerges in a wide variety of natural phenomena and applied sciences including continuum mechanics, thermodynamics, chemical engineering, solid-state electronics, semi-conductor devices, battery research, and nanotechnology. A suitable implementation of a modern methodology, known as Unified Transform Method, is in force throughout this study, with particular challenges arising due to the higher-order mixed-derivative term of the PDE and the generality of the problems under consideration altogether. The boundary and interface conditions appear to be non-standard but are naturally dictated by the structure of the PDE itself. Our explicit analytical formulae directly lend themselves to future explorations of the solutions' qualitative properties such as asymptotic behavior, spatio-temporal dynamics, regularity and well-posedness. This work is expected to be of utility also in the investigation of nonlinear counterparts as well as towards the study of phase-transition phenomena and free-boundary problems, where the interface evolves dynamically according to energy balance laws.

## 1. Introduction

The pseudoparabolic Sobolev-Galpern type [1, 2] equation

(1.1) $$\partial_t u - \alpha \partial_{xxt} u - \beta \partial_{xx} u = 0,$$

which generalizes the diffusion and heat equations, has been proposed as a higher-order model for miscellaneous phenomena and natural processes. Below, for the sake of unity of presentation, we use notations other than those used in the cited papers.

Eq. (1.1) first appeared in 1926 in the paper [3], where monochromatic radiation passing through a gas was considered. The gas contained atoms capable of absorbing radiation. The equation governing the rate of decay of the concentration of excited atoms (atoms that absorbed radiation) had the form of Eq. (1.1). The equation was obtained as a consequence of a system consisting of the balance equation (conservation law)

(1.2) $$\partial_t u + \partial_x v = 0$$

and the constitutive equation

(1.3) $$v = -\alpha \partial_{xt} u - \beta \partial_x u,$$

where $u$ is the concentration of excited atoms and $v$ is the flux of radiation.

In 1960, the same equation (but in three-dimensional space) was derived in the paper [4], but as a model of a completely different phenomenon. The model described seepage (infiltration) of a liquid in fissured rocks. The eqution was obtaind as a consequence of the system of the equations (written here in one-dimensional space)

---

* corresponding author, e-mail address: chatziafrati@math.uoa.gr

$$\alpha\partial_{xx}u + v - u = 0\,,$$

$$\alpha\partial_t vu + \beta(v - u) = 0\,,$$

where $u$ is the average pressure of the liquid in the fissures and $v$ is the average pressure of the liquid in the pores.

Eq. (1.1) was also derived when constructing models of non-Newtonian fluids. In the same period (1960), in the paper [5], the theory of incompressible second-order fluids was developed, and it was shown that several special flow problems lead to Eq. (1.1). The theory was used in the paper [6] to study certain non-steady flows for second-order fluids. It was shown that initial-boundary value problems formulated by physical reasoning turned out to be mathematically well behaved. The Stokes-Rayleigh problem (Stokes's first problem) for Eq. (1.1) was considered in the papers [7, 8]. In the paper [7], the first three terms of the series solution were obtained, while the exact solution to the problem was obtained by Puri in the paper [8]. Several authors criticized Puri's solution as erroneous. However, Christov showed in the paper [9] that Puri's solution is correct, while the critics solutions to this problem are wrong.

Eq. (1.1) was also obtained in theories of heat conduction as one of the generalizations of Fourier's law. In the paper [10], Chen and Gurtin developed a thermodynamic two-temperature theory of heat conduction for the so-called non-simple materials. The linearized energy equation (in one-dimensional space) has the form of the nonhomogeneous equation (1.1). In the paper [11], several generalizations of Fourier's law were proposed in the framework of non-equilibrium thermodynamics with internal variables. The linearized form of one of the constitutive equations is given by

$$v - \alpha\partial_{xx}v = -\beta\partial_x u\,, \tag{1.4}$$

where $u$ is the temperature and $v$ is the heat flux. The energy balance equation (conservation law) is given by Eq. (1.2). Using Eq. (1.2), we conclude that the constitutive equation takes the form of Eq. (1.3). It is necessary to note here that the energy balance equation generally has the form

$$\partial_t u + \partial_x v = f$$

where $f$ is the heat source. In this case, the constitutive equation (1.4) leads to the equation

$$v = -\alpha\partial_{xt}u - \beta\partial_x u + \alpha\partial_x f$$

instead of Eq. (1.3), and the temperature $u$ satisfies the equation

$$\partial_t u - \alpha\partial_{xxt}u - \beta\,\partial_{xx}u = f - \alpha\partial_x f\,,$$

instead of Eq. (1.1).

In more recent advances, the Barenblatt-Zheltov-Kochina model and its variants have arisen in research pertaining to solid-state physics, modern semi-conductor electronics and next-generation lithium-ion rechargeable battery devices; see for instance [12-15] and literature cited therein. More specifically, in advanced battery modelling, pseudo-parabolic or modified diffusion equations appear naturally in several contexts involving complex transport mechanisms in solid-state and liquid-electrolyte systems, especially when models incorporate non-local transport, finite-speed diffusion, or surface-reaction coupling. These equations arise from gradient-energy (capillarity) terms or nonlocal chemical-potential relaxation, which originate in phase-field and Cahn–Hilliard–reaction formulations of intercalation, and they accommodate effects such as stress gradients, microstructural heterogeneity, and finite interface kinetics that cannot be captured by classical Fickian diffusion. In such settings, the effective evolution law often reduces to a Barenblatt-type pseudo-parabolic equation, which is mathematically equivalent to introducing a small-scale spatial memory or inertia in diffusion—an effect that is physically realistic for heterogeneous battery materials, where transport pathways and interfacial processes occur over multiple coupled length scales. Pseudo-parabolic structures also emerge in thermal and electrolyte models when finite-speed diffusion of heat or ions is considered, for example through Maxwell–Cattaneo or non-Fourier heat conduction laws, which are known to be relevant in heterogeneous media like batteries. Taken together, these generalized transport equations provide a unified mathematical framework for incorporating gradient effects, flux relaxation, and nonlocal thermodynamic responses, thereby extending the classical battery models to regimes where sharp concentration fronts, anomalous impedance behavior, or rapid transients play a dominant role. Further research along these lines is currently in progress and will appear in subsequent publications.

A substantial literature exists on the well-posedness, long-time behavior, and other analytical features of Barenblatt's equation and related Sobolev–Galpern type models; see e.g. [16-20] and the extensive bibliography cited therein. Interface problems for linear equations with polynomial dispersion relations have been studied in, e.g., [21,22]. Nevertheless, explicit analytical solutions for interface problems specifically tied to Barenblatt's equation have received no direct attention. This work is essentially intended to address that gap, via suitable implementation of the well-established unified transform method (UTM) pioneered by A.S. Fokas in the late 1990's [23] (see also the works [24-30], the surveys [31,32] and the textbook [33]; for a recent line of applications and further bibliography in connection with this modern PDE technique, see e.g. [34-42], in alphabetical order).

For the benefit and convenience of the reader, we should place the UTM, which is a central tool in this study, into historical perspective. In the second half of the 20th century it was realized that certain nonlinear evolution PDEs, called integrable, can be formulated as the compatibility condition of two linear eigenvalue equations called a Lax pair. This formulation gives rise to a method for solving the initial value problem for these equations, called the inverse scattering transform method. Actually, this method is based on a deeper form of separation of variables. The new approach relied on the panoply of *Riemann-Hilbert* and *d-bar* problems, in synergy with the *Lax pairs* formalism. The UTM could also be applied to certain *non-separable* and *non-self-adjoint* linear problems, as opposed to classical methods such as the ones using the *Laplace*, *Mellin* or (*co*)*sine* transforms. Soon afterwards, an extension was developed which offered a convenient and efficient framework for studying linear evolution and elliptic PDEs with derivatives of arbitrary order in all sorts of spatiotemporal domains. For initial-boundary-value problems (IBVP) for evolution PDE posed on semi-unbounded *half*-spaces and bounded domains, the novel method yields solution representations expressed as contour integrals in the spectral complex plane, much like the traditional *Fourier* transform is involved in the solution of the corresponding counterparts on the *whole* Euclidean space.

In the past 30 years, thanks to a large network of mathematicians, there has been a vast amount of literature on the latter extension of the method to approaching initial-boundary value problems. In fact, the UTM generalises the 'classical' theory, which essentially reduces initial-value problems to Riemann-Hilbert problems via the *scattering transform* and then studies these problems by the method of *inverse scattering* (going back to the works of Gelfand, Levitan and Marchenko). Instead of the usual scattering transform, the new transform is based on the simultaneous spectral analysis of both the $x$-problem and the $t$-problem in the Lax pair. Thus, initial-boundary-value problems are also reduced to Riemann-Hilbert factorization problems in the complex plane. One of the main advantages of the Riemann-Hilbert formulation is that one can use the powerful nonlinear stationary phase and steepest descent theory which gives rigorous results on the asymptotic behavior of solutions to these Riemann-Hilbert problems (as some parameter goes to infinity) and hence it extracts asymptotics for the solution of the associated soliton equation by use of the 'Riemann-Hilbert method'.

The UTM has been successful, both analytically and numerically, in a wide class of linear-PDE problems – of any order and in a variety of domains. At a methodological and pedagogical level, for linear evolution PDEs with polynomial dispersion relation, the UTM typically involves three steps: (1) Derive formally an integral representation for the solution, in the complex Fourier plane. This representation, just like the analogous one obtained via Green's integral identities, contains certain transforms of unknown boundary values; thus it is not yet effective. (2) Analyze the so-called global relation, which is an equation coupling the given initial and boundary data with the unknown boundary values. (3) Use certain invariant properties of the global relation and simple algebraic manipulations to eliminate, from the integral representation obtained in step one, the transforms of the unknown boundary values thereby constructing a candidate solution to the given IBVP.

In recent years, a novel approach to rigorous aspects of the UTM in classical settings has been developed by one of the authors and several collaborators, e.g. [34-41]. That approach has led to valuable new knowledge for a variety of classical IBVPs and PDEs. In addition, it has provided a desirable refinement (extension) and motivated new areas of applicability of the UTM, including the analysis of 'distant' sensitivity (discovery of a long-range instability phenomenon) of solutions of dispersive PDEs, well-posedness and qualitative theory for linear PDEs on the quarter-plane, and so forth. It is here worthwhile noting also that linear PDE have traditionally been of great significance: as abstract objects of theoretical interest, as practical and reasonably accurate models of physical and engineering processes, and as intermediate stages towards investigation of nonlinear problems (in view of e.g. linearizability of certain classes of PDE, travelling-waves analysis, linear estimates necessary for addressing well-posedness questions in low-regularity settings, and so on).

As already mentioned, in this paper, we focus on studying interface problems for Barenblatt's equation on the real line, assuming the interface is located at the origin. Such configurations naturally arise in the study of heat transfer or mass transport through media composed of distinct layers, where physical quantities and structural parameters may undergo abrupt changes. At the point of discontinuity, general interface conditions—determined by the underlying physical model—must be imposed.

Our main objective is to derive explicit solution formulas for these interface problems, thereby offering a fully analytical account of the dynamics in the presence of a jump discontinuity. By applying a suitable extension of the UTM and carrying out a detailed analysis of the resulting expressions, we obtain closed-form integral representations for the solution on both sides of the interface, expressed in terms of the given initial data and interface conditions.

Beyond their intrinsic mathematical value, these explicit solutions provide a platform for future study of free-boundary problems, spatially varying generalizations, nonlinear interface dynamics, numerical methods, and asymptotic behavior in related systems. The techniques introduced here can be adapted to other equations and higher-dimensional settings, offering a versatile analytical approach for problems involving transport across interfaces.

Moreover, our analysis identified a particular class of interface conditions that renders the problem separable: the real-line setting effectively splits into two uncoupled initial–boundary value problems on half-lines. We solve that type of nonstandard half-line IBVP rigorously—examining also its boundary behavior—and subsequently extend the analysis to an analogous finite-interval IBVP. It is noted that more traditional IBVP for Barenblatt's equation with typical boundary conditions (Dirichlet, Neumann, Robin) have been solved via the UTM in [17].

More specifically, we have the following:

***Problem 1*** *Given* $u_0(x)$, $\zeta(t)$ *and* $f(x,t)$, *solve the differential equation*

$$\partial_t u - \alpha\partial_{xxt}u - \beta\,\partial_{xx}u = f\,, \tag{1.5}$$

*for* $u = u(x,t)$, *with* $x > 0$, $t > 0$, *subject to the initial and boundary conditions:*

$$u(x,0) = u_0(x),\ \ x \in \mathbb{R}^+, \tag{1.6}$$

$$\gamma_0[\alpha u_t(0,t) + \beta u(0,t)] + \gamma_1[\alpha u_{xt}(0,t) + \beta u_x(0,t)] = \zeta(t)\,,\ \ t \in \mathbb{R}^+. \tag{1.7}$$

***Assumptions*** We will work with $\alpha > 0$, $\beta > 0$. Also we assume that

$$u_0(x) \in \mathcal{S}([0,\infty)),\ \zeta(t) \in C^\infty([0,\infty)),\ f(x,t) \in C^\infty\left(\overline{\mathbb{R}^+ \times \mathbb{R}^+}\right)\ \textit{and}\ \frac{\partial^m f}{\partial t^m}(\cdot,t) \in \mathcal{S}([0,\infty)),\ \forall m. \tag{1.8}$$

The condition imposed on $f$ means that, for every $T > 0$,

$$\sup\left\{x^\ell\left|\frac{\partial^{n+m} f(x,t)}{\partial x^n \partial t^m}\right| :\ x > 0,\ 0 \le t \le T\right\} < \infty\text{, for all nonnegative integers } n,\ m \text{ and } \ell.$$

The boundary condition (1.7) is – in a sense – natural due to the form of the global relation of equation (1.5) (see (2.1), below) in which there appear the quantities

$$\frac{\lambda e^{-\omega(\lambda)t}}{1+\alpha\lambda^2}\int_{\tau=0}^{t} e^{\lambda\tau}[\alpha u_t(0,\tau) + \beta u(0,\tau)]d\tau\,,\ \ \frac{e^{-\omega(\lambda)t}}{1+\alpha\lambda^2}\int_{\tau=0}^{t} e^{\lambda\tau}[\alpha u_{xt}(0,\tau) + \beta u_x(0,\tau)]d\tau\,.$$

The solution of Problem 1 and its properties near the axes are established in Sections 2 and 3.

We were led to the above problem by considering the following interface problem:

***Problem 2*** *Solve the initial value problems*

$$\begin{cases} \partial_t u^R - \alpha_R\partial_{xxt}u^R - \beta_R\,\partial_{xx}u^R = f_R, \\ u^R(x,0) = u_0^R(x),\ \ x \in \mathbb{R}^+, \end{cases} \tag{1.9}$$

*for* $u^R = u^R(x,t)$, *with* $x > 0$, $t > 0$, *and*

$$\begin{cases} \partial_t u^L - \alpha_L\partial_{xxt}u^L - \beta_L\,\partial_{xx}u^L = f_L, \\ u^L(x,0) = u_0^L(x),\ \ x \in \mathbb{R}^-, \end{cases} \tag{1.10}$$

*for* $u^L = u^L(x,t)$, *with* $x < 0$, $t > 0$, *subject to the following interface conditions: For* $t > 0$,

(1.11) $\gamma_{11}[\alpha_R u_t^R(0,t)+\beta_R u^R(0,t)]+\gamma_{12}[\alpha_R u_{xt}^R(0,t)+\beta_R u_x^R(0,t)]$
$$+\gamma_{13}[\alpha_L u_t^L(0,t)+\beta_L u^L(0,t)]+\gamma_{14}[\alpha_L u_{xt}^L(0,t)+\beta_L u_x^L(0,t)]=I_1(t),$$

(1.12) $\gamma_{21}[\alpha_R u_t^R(0,t)+\beta_R u^R(0,t)]+\gamma_{22}[\alpha_R u_{xt}^R(0,t)+\beta_R u_x^R(0,t)]$
$$+\gamma_{23}[\alpha_L u_t^L(0,t)+\beta_L u^L(0,t)]+\gamma_{24}[\alpha_L u_{xt}^L(0,t)+\beta_L u_x^L(0,t)]=I_2(t).$$

We will prove in section 4 that problem 2 is *«separable»* in the sense that in order for this problem to be solvable we have to assume that

$$\Gamma_{12}:=\det\begin{bmatrix}\gamma_{11} & \gamma_{12}\\ \gamma_{21} & \gamma_{22}\end{bmatrix}=0 \text{ and } \Gamma_{34}:=\det\begin{bmatrix}\gamma_{13} & \gamma_{14}\\ \gamma_{23} & \gamma_{24}\end{bmatrix}=0. \tag{1.13}$$

It is easy to see that assumption (1.9) imply that problem 2 is reduced essentially to this problem with the interface conditions simplified to the following form:

$$\gamma_{21}[\alpha_R u_t^R(0,t)+\beta_R u^R(0,t)]+\gamma_{22}[\alpha_R u_{xt}^R(0,t)+\beta_R u_x^R(0,t)]=I_1(t)$$

$$\gamma_{13}[\alpha_L u_t^L(0,t)+\beta_L u^L(0,t)]+\gamma_{14}[\alpha_L u_{xt}^L(0,t)+\beta_L u_x^L(0,t)]=I_2(t),$$

i.e., to two separate IBV problems, one in the right half-line and another one in the left half-line.

In section 2 we implement the UTM in order to derive the solution of problem 1 (see formula (1.14), below) and in section 3 we verify that indeed this is the solution.

*Note* Setting $v(x,t)=\gamma_0 u(x,t)+\gamma_1 u_x(x,t)$, we see that Problem 1 is reduced to the following problem

$$\begin{cases}\partial_t v-\alpha\partial_{xxt}v-\beta\partial_{xx}v=f,\ x>0,\ t>0,\\ v(x,0)=\gamma_0 u_0(x)+\gamma_1 u_0'(x),\ x>0,\\ \alpha v_t(0,t)+\beta v(0,t)=\zeta(t),\ t>0,\end{cases}$$

for $v(x,t)$. The above problem is the special case of problem 1, corresponding to $\gamma_1=0$. If we chose to solve the above problem, after finding $v(x,t)$, we would have to solve, in addition, the differential equations

$$\gamma_0 u(x,t)+\gamma_1 u_x(x,t)=v(x,t), \text{ for each fixed } t>0.$$

We considered this approach more complicated, especially as far as the verification process is concerned, and we preferred, instead, to solve Problem 1 directly in the general case, i.e., with $\gamma_1\neq 0$.

Section 5 deals with the derivation of the UTM solution for very general interface problems on the line. (For the precise statement, see problem 3, described in Section 5.)

Finally, in Section 6, we consider an initial boundary value problem for this equation on a finite interval $\{x:\ 0<x<l\}$. More precisely we derive the UTM solution to the following problem:

***Problem 4*** *Given $u_0(x)$, $\zeta_0(t)$, $\zeta_l(t)$ and $f(x,t)$, solve the differential equation*

$$\partial_t u-\alpha\partial_{xxt}u-\beta\partial_{xx}u=f,$$

*for $u=u(x,t)$, with $0<x<l$, $t>0$, subject to the initial and boundary conditions:*

$$u(x,0)=u_0(x),\ \ 0<x<l, \tag{1.14}$$

*and, for $t\in\mathbb{R}^+$,*

$$\begin{cases}[\alpha u_{xt}(0,t)+\beta u_x(0,t)]-c_0 u(0,t)=\zeta_0(t)\\ [\alpha u_{xt}(l,t)+\beta u_x(l,t)]+c_l u(l,t)=\zeta_l(t).\end{cases} \tag{1.15}$$

The boundary conditions can be considered as an analog of the Robin [43] boundary conditions for the heat equation. Indeed, the equation can be written as

$$\partial_t u + \partial_x v = f(x,t)\,,\ v = -(\alpha \partial_{xt} u + \beta \partial_x u)\,.$$

This equation is in fact the balance equation (conservation law), where $v$ is the flux of $u$ and $f$ is the sourse term. If $c_0 = 0$ and $c_l = 0$, then the boundary conditions become an analog of the Neumann boundary conditions for the heat equation. If $\alpha = 0$, then Problem 4 becomes the initial-boundary value problem with the Robin boundary conditions for the heat equation.

***Notation*** We set

$$\omega = \omega(\lambda) = \frac{\beta\lambda^2}{1+\alpha\lambda^2}$$

so that

$$\left(\frac{\partial}{\partial t} - \alpha\frac{\partial^3}{\partial x^2 \partial t} - \beta\frac{\partial^2}{\partial x^2}\right) e^{i\lambda x - \omega(\lambda)t} = 0\,. \tag{1.16}$$

In order to justify the convergence of certain integrals as well as the deformation of contours, it is important to notice that

$$\lim_{\lambda\to\infty} \omega(\lambda) = \lim_{\lambda\to\infty} \frac{\beta\lambda^2}{1+\alpha\lambda^2} = \frac{\beta}{\alpha}\,.$$

Hence,

$$\sup_{\left|1+\alpha\lambda^2\right|\ge\varepsilon} \left|e^{\pm\omega(\lambda)}\right| = \sup_{\left|1+\alpha\lambda^2\right|\ge\varepsilon} \left|\exp\left(\pm\frac{\beta\lambda^2}{1+\alpha\lambda^2}\right)\right| < +\infty \ \ (\forall \varepsilon > 0\,). \tag{1.17}$$

We will also use the following notation of transforms:

$$\hat{u}_0(\lambda) = \int_0^\infty e^{-i\lambda x} u_0(x) dx\,,\ \lambda\in\mathbb{C}\,,\ \operatorname{Im}\lambda \le 0\,,$$

$$\hat{u}(\lambda,t) = \int_0^\infty e^{-i\lambda x} u(x,t) dx\,,\ \lambda\in\mathbb{C}\,,\ \operatorname{Im}\lambda \le 0\,,$$

$$\tilde{\zeta}(\omega,t) = \int_0^t e^{\omega(\lambda)\tau} \zeta(\tau) d\tau\,,\ \lambda\in\mathbb{C}\,,$$

$$\hat{f}(\lambda,t) = \int_0^\infty e^{-i\lambda x} f(x,t) dx\,,\ \lambda\in\mathbb{C}\,,\ \operatorname{Im}\lambda \le 0\,,$$

$$\tilde{\hat{f}}(\lambda,\omega,t) = \int_0^t e^{\omega(\lambda)\tau} \hat{f}(\lambda,\tau) d\tau\,,\ \lambda\in\mathbb{C}\,,\ \operatorname{Im}\lambda \le 0\,.$$

Finally, setting

$$\mathcal{P} := \gamma_0\left(\alpha\frac{\partial}{\partial t} + \beta\right) + \gamma_1\left(\alpha\frac{\partial^2}{\partial x \partial t} + \beta\frac{\partial}{\partial x}\right),$$

we may write the boundary condition (1.3) as follows: $(\mathcal{P}u)(0,t) = \zeta(t)\,,\ t > 0\,.$

***Solution of problem 1*** We will show that the UTM solution (derived in Section 2) of *problem 1* is given by the following formula: For $x > 0$ and $t > 0$,

$$u(x,t) = \frac{1}{2\pi}\int_{-\infty}^{\infty} \hat{u}_0(\lambda) e^{i\lambda x - \omega(\lambda)t} d\lambda - \frac{1}{2\pi}\int_{\mathcal{K}(i/\sqrt{\alpha})} \hat{u}_0(-\lambda) \frac{\gamma_0 - i\gamma_1\lambda}{\gamma_0 + i\gamma_1\lambda} e^{i\lambda x - \omega(\lambda)t} d\lambda \tag{1.18}$$

$$-\frac{i}{\pi}\int_{\mathcal{K}(i/\sqrt{\alpha})}\frac{e^{i\lambda x-\omega(\lambda)t}}{(1+\alpha\lambda^2)(i\gamma_1\lambda+\gamma_0)}\lambda\tilde{\zeta}(\omega,t)d\lambda$$
$$+\frac{1}{2\pi}\int_{-\infty}^{\infty}\frac{e^{i\lambda x-\omega(\lambda)t}}{1+\alpha\lambda^2}\tilde{\hat{f}}(\lambda,\omega,t)d\lambda-\frac{1}{2\pi}\int_{\mathcal{K}(i/\sqrt{\alpha})}\frac{e^{i\lambda x-\omega(\lambda)t}}{1+\alpha\lambda^2}\frac{\gamma_0-i\gamma_1\lambda}{\gamma_0+i\gamma_1\lambda}\tilde{\hat{f}}(-\lambda,\omega,t)d\lambda,$$

where $\mathcal{K}(i/\sqrt{\alpha})$ is a sufficiently small simple closed curve surrounding the point $i/\sqrt{\alpha}$, traversed counterclockwise (fig. 1).

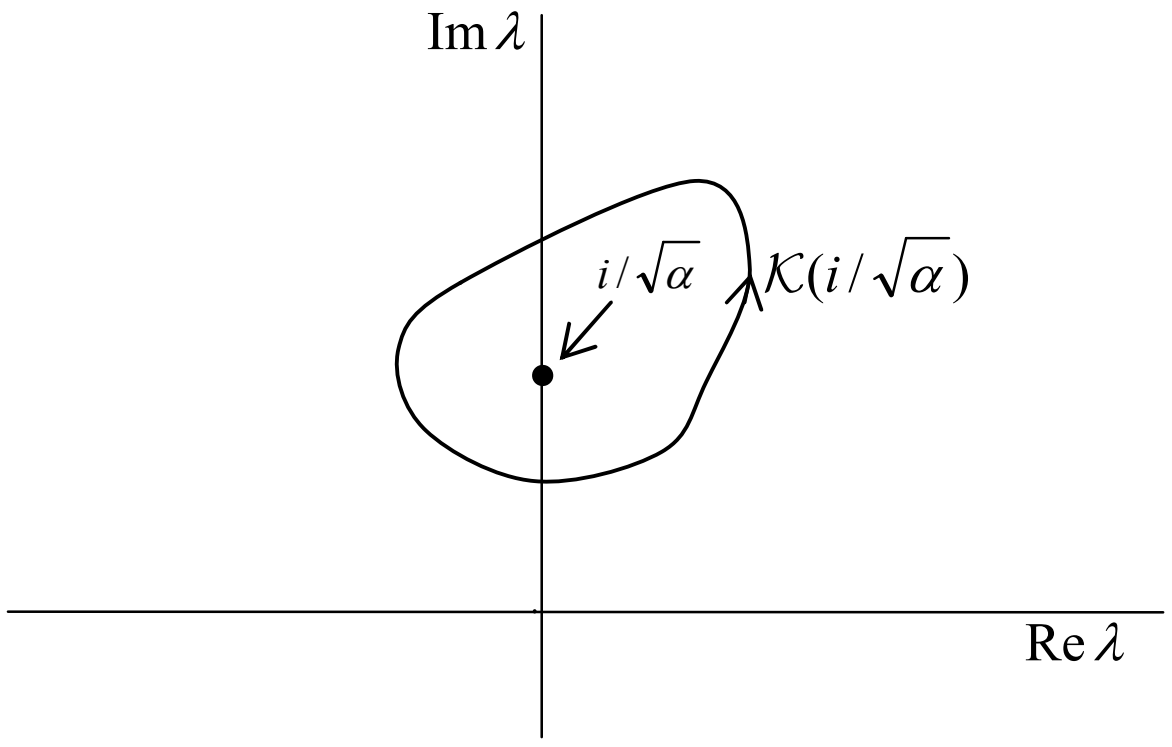

**Fig. 1** A typical contour $\mathcal{K}(i/\sqrt{\alpha})$ and its orientaion

More precisely, for Problem 1 we have:

***Theorem 1*** *Assuming that the data $u_0(x)$, $\zeta(t)$ and $f(x,t)$ satisfy (1.4), the function $u(x,t)$, defined by (1.14), is $C^\infty$ (jointly) for $(x,t)\in\mathbb{R}^+\times\mathbb{R}^+$, satisfies the differential equation (1.1),*

*(1.19)* $$\lim_{t\to0^+}u(x,t)=u_0(x)\text{, uniformly for } x>0,$$

*and, for every $T>0$,*

*(1.20)* $$\lim_{x\to0^+}\{\gamma_0[\alpha u_t(x,t)+\beta u(x,t)]+\gamma_1[\alpha u_{xt}(x,t)+\beta u_x(x,t)]\}=\zeta(t)\text{, uniformly for } 0<t\le T.$$

***Theorem 2*** *The function $u(x,t)$, defined for $(x,t)\in\mathbb{R}^+\times\mathbb{R}^+$ by (1.14), satisfies the following:*

*1st The limits*

$$U_{n,m}(x):=\lim_{t\to0^+}\frac{\partial^{n+m}u(x,t)}{\partial x^n\partial t^m},\ n,m=0,1,2,...,$$

*exist and are uniform for $x>0$. Moreover the functions $U_{n,m}(x)$ are $C^\infty$ for $x\in\mathbb{R}^+$.*

*2nd The limits*

$$H_{n,m}(t):=\lim_{x\to0^+}\frac{\partial^{n+m}u(x,t)}{\partial x^n\partial t^m},\ n,m=0,1,2,...,$$

*exist and are uniform for $0<t\le T$ (for every $T>0$). Moreover the functions $H_{n,m}(t)$ are $C^\infty$ for $t\in\mathbb{R}^+$.*

*3rd The limit* $\lim_{t\to0^+}\frac{\partial^n u(x,t)}{\partial x^n}=\frac{d^n u_0(x)}{dx^n}$*, uniformly for $x>0$.*

*4th The limit*

$$\lim_{x\to0^+}\left\{\gamma_0\left[\alpha\frac{\partial^{m+1}u(x,t)}{\partial t^{m+1}}+\beta\frac{\partial^m u(x,t)}{\partial t^m}\right]+\gamma_1\left[\alpha\frac{\partial^{m+2}u(x,t)}{\partial t^{m+1}\partial x}+\beta\frac{\partial^{m+1}u(x,t)}{\partial t^m\partial x}\right]\right\}=\frac{d^m\zeta(t)}{dt^m}\text{, uniformly for } 0<t\le T.$$

***Theorem 3*** *For* $T>0$, *the function* $u(x,t)$, *defined for* $(x,t)\in\mathbb{R}^+\times\mathbb{R}^+$ *by (1.14), satisfies,*

$$\lim_{x\to+\infty}\left(x^{\ell}\frac{\partial^{n+m}u(x,t)}{\partial x^n\partial t^m}\right)=0,$$

*for nonnegative integers* $n$, $m$ *and* $\ell$, *uniformly for* $0<t\le T$.

The proofs of the above theorems are discussed in Section 3.

## 2. Derivation of the solution of Problem 1 via the UTM

Working as in [17], we derive the following equation, known as the *global relation*:

$$(2.1)\quad \hat{u}(\lambda,t)=\hat{u}_0(\lambda)e^{-\omega(\lambda)t}-\frac{i\lambda e^{-\omega(\lambda)t}}{1+\alpha\lambda^2}[\alpha\tilde{g}'_0(\omega,t)+\beta\tilde{g}_0(\omega,t)]$$
$$-\frac{e^{-\omega(\lambda)t}}{1+\alpha\lambda^2}[\alpha\tilde{g}'_1(\omega,t)+\beta\tilde{g}_1(\omega,t)]+\frac{e^{-\omega(\lambda)t}}{1+\alpha\lambda^2}\tilde{\hat{f}}(\lambda,\omega,t),\ \operatorname{Im}\lambda\le 0,$$

where $g_0(t):=u(0,t)$, $g'_0(t):=u_t(0,t)$, $g_1(t):=u_x(0,t)$ and $g'_1(t):=u_{xt}(0,t)$.

Integrating (2.1) we obtain

$$u(x,t)=\frac{1}{2\pi}\int_{-\infty}^{\infty}e^{i\lambda x}\hat{u}(\lambda,t)d\lambda$$
$$=\frac{1}{2\pi}\int_{-\infty}^{\infty}\hat{u}_0(\lambda)e^{i\lambda x-\omega(\lambda)t}d\lambda-\frac{1}{2\pi}\int_{-\infty}^{\infty}\frac{i\lambda e^{i\lambda x-\omega(\lambda)t}}{1+\alpha\lambda^2}[\alpha\tilde{g}'_0(\omega,t)+\beta\tilde{g}_0(\omega,t)]d\lambda$$
$$-\frac{1}{2\pi}\int_{-\infty}^{\infty}\frac{e^{i\lambda x-\omega(\lambda)t}}{1+\alpha\lambda^2}[\alpha\tilde{g}'_1(\omega,t)+\beta\tilde{g}_1(\omega,t)]d\lambda+\frac{1}{2\pi}\int_{-\infty}^{\infty}\frac{e^{i\lambda x-\omega(\lambda)t}}{1+\alpha\lambda^2}\tilde{\hat{f}}(\lambda,\omega,t)d\lambda,\ x>0,\ t>0.$$

In view of (1.17), the integrals in the above equation, which involve the terms $\alpha\tilde{g}'_0(\omega,t)+\beta\tilde{g}_0(\omega,t)$ and $\alpha\tilde{g}'_1(\omega,t)+\beta\tilde{g}_1(\omega,t)$, can be transformed, leading to

$$(2.2)\quad u(x,t)=\frac{1}{2\pi}\int_{-\infty}^{\infty}\hat{u}_0(\lambda)e^{i\lambda x-\omega(\lambda)t}d\lambda-\frac{1}{2\pi}\int_{\lambda\in\mathcal{K}(i/\sqrt{\alpha})}\frac{i\lambda e^{i\lambda x-\omega(\lambda)t}}{1+\alpha\lambda^2}[\alpha\tilde{g}'_0(\omega,t)+\beta\tilde{g}_0(\omega,t)]d\lambda$$
$$-\frac{1}{2\pi}\int_{\lambda\in\mathcal{K}(i/\sqrt{\alpha})}\frac{e^{i\lambda x-\omega(\lambda)t}}{1+\alpha\lambda^2}[\alpha\tilde{g}'_1(\omega,t)+\beta\tilde{g}_1(\omega,t)]d\lambda+\frac{1}{2\pi}\int_{-\infty}^{\infty}\frac{e^{i\lambda x-\omega(\lambda)t}}{1+\alpha\lambda^2}\tilde{\hat{f}}(\lambda,\omega,t)d\lambda,\ x>0,\ t>0.$$

Also, (2.1) gives

$$\hat{u}(-\lambda,t)=\hat{u}_0(-\lambda)e^{-\omega(\lambda)t}+\frac{i\lambda e^{-\omega(\lambda)t}}{1+\alpha\lambda^2}[\alpha\tilde{g}'_0(\omega,t)+\beta\tilde{g}_0(\omega,t)]$$
$$-\frac{e^{-\omega(\lambda)t}}{1+\alpha\lambda^2}[\alpha\tilde{g}'_1(\omega,t)+\beta\tilde{g}_1(\omega,t)]+\frac{e^{-\omega(\lambda)t}}{1+\alpha\lambda^2}\tilde{\hat{f}}(-\lambda,\omega,t),\ \operatorname{Im}\lambda\ge 0.$$

Equivalently,

$$(2.3)\quad (1+\alpha\lambda^2)e^{\omega(\lambda)t}\hat{u}(-\lambda,t)=(1+\alpha\lambda^2)\hat{u}_0(-\lambda)+i\lambda[\alpha\tilde{g}'_0(\omega,t)+\beta\tilde{g}_0(\omega,t)]$$
$$-[\alpha\tilde{g}'_1(\omega,t)+\beta\tilde{g}_1(\omega,t)]+\tilde{\hat{f}}(-\lambda,\omega,t),\ \operatorname{Im}\lambda\ge 0.$$

Setting

$G_0^*(\lambda,t):=\alpha\tilde{g}'_0(\omega,t)+\beta\tilde{g}_0(\omega,t)$, $G_1^*(\lambda,t):=\alpha\tilde{g}'_1(\omega,t)+\beta\tilde{g}_1(\omega,t)$,

$J(\lambda,t):=-(1+\alpha\lambda^2)\hat{u}_0(-\lambda)-\tilde{\hat{f}}(-\lambda,\omega,t)$ and $J^*(\lambda,t):=(1+\alpha\lambda^2)e^{\omega(\lambda)t}\hat{u}(-\lambda,t)$,

we write (2.3) as follows:

$$i\lambda G_0^*(\lambda,t)-G_1^*(\lambda,t)=J(\lambda,t)+J^*(\lambda,t)\,. \tag{2.4}$$

Then (2.2) becomes

$$u(x,t)=\frac{1}{2\pi}\int_{-\infty}^{\infty}\hat{u}_0(\lambda)e^{i\lambda x-\omega(\lambda)t}d\lambda-\frac{1}{2\pi}\int_{\lambda\in\mathcal{K}(i/\sqrt{\alpha})}\frac{i\lambda e^{i\lambda x-\omega(\lambda)t}}{1+\alpha\lambda^2}G_0^*(\lambda,t)d\lambda$$
$$-\frac{1}{2\pi}\int_{\lambda\in\mathcal{K}(i/\sqrt{\alpha})}\frac{e^{i\lambda x-\omega(\lambda)t}}{1+\alpha\lambda^2}G_1^*(\lambda,t)d\lambda+\frac{1}{2\pi}\int_{-\infty}^{\infty}\frac{e^{i\lambda x-\omega(\lambda)t}}{1+\alpha\lambda^2}\tilde{\hat{f}}(\lambda,\omega,t)d\lambda\,,\ x>0\,,\ t>0\,. \tag{2.5}$$

On the other hand, (1.7) gives the equation

$$\gamma_0 G_0^*(\lambda,t)+\gamma_1 G_1^*(\lambda,t)=\tilde{\zeta}(\omega,t)\,. \tag{2.6}$$

Solving the system of equations (2.4) & (2.6), we find

$$G_0^*(\lambda,t)=\frac{\mathfrak{D}_0^*(\lambda,t)}{\mathfrak{D}(\lambda)}\ \text{and}\ G_1^*(\lambda,t)=\frac{\mathfrak{D}_1^*(\lambda,t)}{\mathfrak{D}(\lambda)}\,,$$

where

$$\mathfrak{D}(\lambda):=\det\begin{bmatrix} i\lambda & -1\\ \gamma_0 & \gamma_1\end{bmatrix}=i\gamma_1\lambda+\gamma_0\,,$$

$$\mathfrak{D}_0^*(\lambda,t):=\det\begin{bmatrix} J(\lambda,t)+J^*(\lambda,t) & -1\\ \tilde{\zeta}(\omega,t) & \gamma_1\end{bmatrix}=\gamma_1[J(\lambda,t)+J^*(\lambda,t)]+\tilde{\zeta}(\omega,t)\,,$$

$$\mathfrak{D}_1^*(\lambda,t):=\det\begin{bmatrix} i\lambda & J(\lambda,t)+J^*(\lambda,t)\\ \gamma_0 & \tilde{\zeta}(\omega,t)\end{bmatrix}=-\gamma_0[J(\lambda,t)+J^*(\lambda,t)]+i\lambda\tilde{\zeta}(\omega,t)\,.$$

Assuming that $\mathfrak{D}(i/\sqrt{\alpha})\neq 0$, i.e., $\gamma_0\sqrt{\alpha}\neq\gamma_1$, we claim that

$$\int_{\lambda\in\mathcal{K}(i/\sqrt{\alpha})} i\lambda e^{i\lambda x-\omega(\lambda)t}\frac{J^*(\lambda,t)}{(1+\alpha\lambda^2)(i\gamma_1\lambda+\gamma_0)}d\lambda=0\ \text{ and }\ \int_{\lambda\in\mathcal{K}(i/\sqrt{\alpha})} e^{i\lambda x-\omega(\lambda)t}\frac{J^*(\lambda,t)}{(1+\alpha\lambda^2)(i\gamma_1\lambda+\gamma_0)}d\lambda=0\,, \tag{2.7}$$

provided that the contour $\mathcal{K}(i/\sqrt{\alpha})$, in the case $\gamma_1\neq 0$, is chosen so that

$$\frac{i\gamma_0}{\gamma_1}\notin \operatorname{int}[\mathcal{K}(i/\sqrt{\alpha})]\,.$$

Indeed, (2.7) follows from Cauchy's theorem since

$$e^{i\lambda x-\omega(\lambda)t}\frac{J^*(\lambda,t)}{(1+\alpha\lambda^2)(i\gamma_1\lambda+\gamma_0)}=e^{i\lambda x}\frac{\hat{u}(-\lambda,t)}{i\gamma_1\lambda+\gamma_0}\,.$$

Thus, in view of (2.7), formula (2.5) gives:

$$u(x,t)=\frac{1}{2\pi}\int_{-\infty}^{\infty}\hat{u}_0(\lambda)e^{i\lambda x-\omega(\lambda)t}d\lambda-\frac{1}{2\pi}\int_{\lambda\in\mathcal{K}(i/\sqrt{\alpha})} i\lambda e^{i\lambda x-\omega(\lambda)t}\frac{\gamma_1 J(\lambda,t)+\tilde{\zeta}(\omega,t)}{(1+\alpha\lambda^2)(i\gamma_1\lambda+\gamma_0)}d\lambda$$
$$-\frac{1}{2\pi}\int_{\lambda\in\mathcal{K}(i/\sqrt{\alpha})} e^{i\lambda x-\omega(\lambda)t}\frac{-\gamma_0 J(\lambda,t)+i\lambda\tilde{\zeta}(\omega,t)}{(1+\alpha\lambda^2)(i\gamma_1\lambda+\gamma_0)}d\lambda+\frac{1}{2\pi}\int_{-\infty}^{\infty}\frac{e^{i\lambda x-\omega(\lambda)t}}{1+\alpha\lambda^2}\tilde{\hat{f}}(\lambda,\omega,t)d\lambda\,,\ x>0\,,\ t>0\,.$$

Substituting $J(\lambda,t)$ by its definition, we see that the above equation becomes (1.18).

## 3. Proof of Theorems 1, 2 and 3

We split the proof of Theorem 1 in several steps.

***Step 1*** We claim that the integrals in (1.18) define $C^\infty$ functions for $(x,t)\in\mathbb{R}^+\times\mathbb{R}^+$. This is immediate for the integrals taken on the contour $\mathcal{K}(i/\sqrt{\alpha})$. However the integral

(3.1) $$\int_{-\infty}^{\infty}\hat{u}_0(\lambda)e^{i\lambda x-\omega(\lambda)t}d\lambda,$$

in general, does not even converge absolutely, and becomes worse if we differentiate it with respect to $x$. To deal with this integral we use the identity

(3.2) $$\hat{u}_0(\lambda)=\sigma_M(\lambda)+\frac{1}{(i\lambda)^M}\int_0^\infty e^{-i\lambda y}u_0^{(M)}(y)dy \text{ for } \lambda\in\mathbb{C}-\{0\} \text{ with } \operatorname{Im}\lambda\le 0,$$

where we have set

$$\sigma_M(\lambda)=\sum_{k=1}^{M}\frac{u_0^{(k-1)}(0)}{(i\lambda)^k}.$$

Let us observe that $\hat{u}_0(\lambda)-\sigma_M(\lambda)=\mathrm{O}(1/\lambda^{M+1})$, as $\lambda\to\infty$ with $\operatorname{Im}\lambda\le 0$.

Thus, for every $M$, the integral (3.1) can be written as follows:

(3.3) $$\int_{-1}^{1}e^{i\lambda x-\omega(\lambda)t}\hat{u}_0(\lambda)d\lambda+\sum_{k=1}^{M}u_0^{(k-1)}(0)\left(\int_{-\infty}^{-1}+\int_{1}^{\infty}\right)e^{i\lambda x-\omega(\lambda)t}\frac{d\lambda}{(i\lambda)^k}+\left(\int_{-\infty}^{-1}+\int_{1}^{\infty}\right)e^{i\lambda x-\omega(\lambda)t}(u_0^{(M)})\hat{}(\lambda)\frac{d\lambda}{(i\lambda)^M}.$$

But in view of (1.17), Cauchy's theorem and Jordan's lemma imply that

(3.4) $$\left(\int_{-\infty}^{-1}+\int_{1}^{\infty}\right)e^{i\lambda x-\omega(\lambda)t}\frac{d\lambda}{(i\lambda)^k}=\int_{\mathcal{C}}e^{i\lambda x-\omega(\lambda)t}\frac{d\lambda}{(i\lambda)^k}\quad(k\ge 1),$$

where $\mathcal{C}$ is any simple contour, starting at 1, ending at $-1$, lying (except its endpoints) in the domain $\{\operatorname{Im}\lambda>0,|\lambda|>1/\sqrt{\alpha}\}$, and passing "above" the point $i/\sqrt{\alpha}$. (See fig 2.)

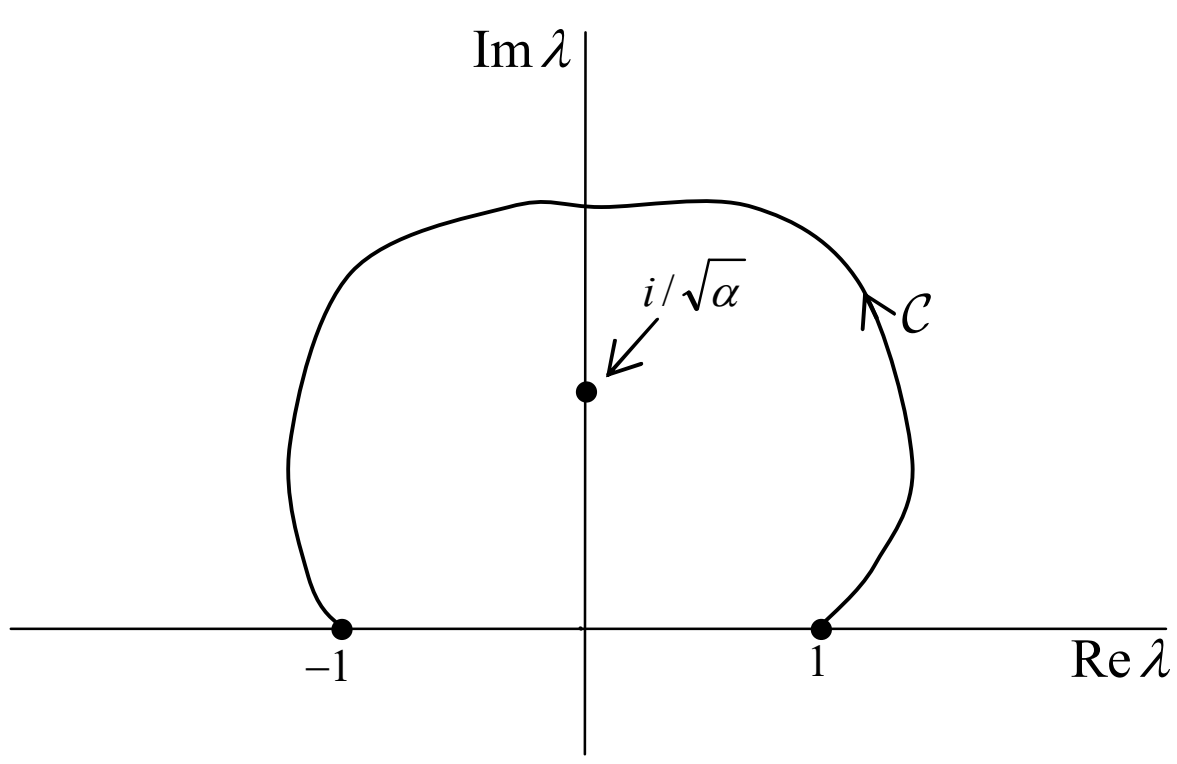

**Fig. 2** A typical contour $\mathcal{C}$

Since (3.3) holds for every $M$, substituting (3.4) in (3.3), we easily obtain that the integral (3.1) defines a $C^\infty$ function in $(x,t)\in\mathbb{R}^+\times\mathbb{R}^+$. Indeed, given $n$ and $m$, we have

(3.5) $$\frac{\partial^{n+m}}{\partial x^n\partial t^m}\left[\int_{\lambda=-\infty}^{\infty}e^{i\lambda x-\omega(\lambda)t}\hat{u}_0(\lambda)d\lambda\right]=\int_{-1}^{1}(i\lambda)^n(-\omega)^m e^{i\lambda x-\omega t}\hat{u}_0(\lambda)d\lambda$$

$$+\sum_{k=1}^{M} u_0^{(k-1)}(0)\int_{\lambda\in\mathcal{C}}(-\omega)^m e^{i\lambda x-\omega t}\frac{d\lambda}{(i\lambda)^{k-n}}+\left(\int_{-\infty}^{-1}+\int_{1}^{\infty}\right)(-\omega)^m e^{i\lambda x-\omega t}(u_0^{(M)}\hat{)}(\lambda)\frac{d\lambda}{(i\lambda)^{M-n}},$$

provided that $M$ is sufficiently large.

To deal with the integral involving the term $\hat{f}(\lambda,t)$, we use the identity

$$\hat{f}(\lambda,t)=\xi_M(\lambda,t)+\frac{1}{(i\lambda)^M}\int_0^\infty e^{-i\lambda y}\frac{\partial^M f(y,t)}{\partial y^M}dy\text{ , where }\xi_M(\lambda,t)=\sum_{k=1}^{M}\frac{1}{(i\lambda)^k}\left.\frac{\partial^{k-1}f(y,t)}{\partial y^{k-1}}\right|_{y=0},$$

which generalizes (3.2). Then

$$\eta_M(\lambda,t):=\hat{f}(\lambda,t)-\xi_M(\lambda,t)=\mathrm{O}(1/\lambda^{M+1})\text{, for }\lambda\to\infty\text{, with }\operatorname{Im}\lambda\le 0, \tag{3.6}$$

uniformly for $t$ in compact subsets of $[0,+\infty)$.

Thus

$$\tilde{\hat{f}}(\lambda,\omega(\lambda),t)=\int_0^t e^{\omega(\lambda)\tau}\eta_M(\lambda,\tau)d\tau+\int_0^t e^{\omega(\lambda)\tau}\xi_M(\lambda,\tau)d\tau=\tilde{\eta}_M(\lambda,\omega(\lambda),t)+\tilde{\xi}_M(\lambda,\omega(\lambda),t). \tag{3.7}$$

Although (3.7) holds for $\lambda\in\mathbb{C}-\{0\}$ with $\operatorname{Im}\lambda\le 0$, the last integral in (3.7), i.e. $\tilde{\xi}_M(\lambda,\omega(\lambda),t)$, is defined for every $\lambda\in\mathbb{C}-\{0\}$.

Therefore, for $x>0$ and $t>0$,

$$\int_{-\infty}^{\infty}e^{i\lambda x-\omega(\lambda)t}\frac{1}{1+\alpha\lambda^2}\tilde{\hat{f}}(\lambda,\omega(\lambda),t)d\lambda=\int_{-1}^{1}e^{i\lambda x-\omega(\lambda)t}\frac{1}{1+\alpha\lambda^2}\tilde{\hat{f}}(\lambda,\omega(\lambda),t)d\lambda$$
$$+\int_{\mathcal{C}}e^{i\lambda x-\omega(\lambda)t}\frac{1}{1+\alpha\lambda^2}\tilde{\xi}_M(\lambda,\omega(\lambda),t)d\lambda+\left(\int_{-\infty}^{-1}+\int_{1}^{\infty}\right)e^{i\lambda x-\omega(\lambda)t}\frac{1}{1+\alpha\lambda^2}\tilde{\eta}_M(\lambda,\omega(\lambda),t)d\lambda. \tag{3.8}$$

It follows that the integral in the LHS of (3.8) is $C^\infty$ in $(x,t)\in\mathbb{R}^+\times\mathbb{R}^+$.

Furthermore, the differentiations $\partial^{n+m}/\partial x^n\partial t^m$ can be interchanged with the integrals in the RHS of (3.8), leading to formulas for the derivatives:

$$\frac{\partial^{n+m}}{\partial x^n\partial t^m}\left[\int_{-\infty}^{\infty}e^{i\lambda x-\omega(\lambda)t}\frac{1}{1+\alpha\lambda^2}\tilde{\hat{f}}(\lambda,\omega(\lambda),t)d\lambda\right]=\int_{-1}^{1}e^{i\lambda x-\omega(\lambda)t}\frac{(i\lambda)^n(-\omega)^m}{1+\alpha\lambda^2}\tilde{\hat{f}}(\lambda,\omega(\lambda),t)d\lambda$$
$$+\int_{\mathcal{C}}e^{i\lambda x-\omega(\lambda)t}\frac{(i\lambda)^n(-\omega)^m}{1+\alpha\lambda^2}\tilde{\xi}_M(\lambda,\omega(\lambda),t)d\lambda+\left(\int_{-\infty}^{-1}+\int_{1}^{\infty}\right)e^{i\lambda x-\omega(\lambda)t}\frac{(i\lambda)^n(-\omega)^m}{1+\alpha\lambda^2}\tilde{\eta}_M(\lambda,\omega(\lambda),t)d\lambda, \tag{3.9}$$

provided that $M$ is sufficiently large.

***Step 2*** We claim that the function $u(x,t)$ satisfies the differential equation (1.5). In view of (1.16) and the computations made in *Step 1*, it suffices to show the following:

$$\int_{\mathcal{K}(i/\sqrt{\alpha})}\left\{\frac{1}{(1+\alpha\lambda^2)(i\gamma_1\lambda+\gamma_0)}\left(e^{i\lambda x-\omega(\lambda)t}\frac{\partial}{\partial t}[\tilde{\zeta}(\omega,t)]-\alpha\frac{\partial}{\partial t}[\tilde{\zeta}(\omega,t)]\frac{\partial^2}{\partial x^2}[e^{i\lambda x-\omega(\lambda)t}]\right)\right\}\lambda d\lambda=0, \tag{3.10}$$

$$\int_{-\infty}^{\infty}\frac{1}{1+\alpha\lambda^2}\left(e^{i\lambda x-\omega(\lambda)t}\frac{\partial}{\partial t}[\tilde{\hat{f}}(\lambda,\omega,t)]-\alpha\frac{\partial}{\partial t}[\tilde{\hat{f}}(\lambda,\omega,t)]\frac{\partial^2}{\partial x^2}[e^{i\lambda x-\omega(\lambda)t}]\right)d\lambda=2\pi f(x,t), \tag{3.11}$$

$$\int_{\mathcal{K}(i/\sqrt{\alpha})}\frac{e^{i\lambda x-\omega(\lambda)t}}{1+\alpha\lambda^2}\frac{\gamma_0-i\gamma_1\lambda}{\gamma_0+i\gamma_1\lambda}\left(\frac{\partial}{\partial t}[\tilde{\hat{f}}(-\lambda,\omega,t)]-\alpha\frac{\partial}{\partial t}[\tilde{\hat{f}}(-\lambda,\omega,t)]\frac{\partial^2}{\partial x^2}[e^{i\lambda x-\omega(\lambda)t}]\right)d\lambda=0. \tag{3.12}$$

*Proof of (3.10)* Since

$$\frac{\partial}{\partial t}[\tilde{\zeta}(\omega,t)]-\alpha\frac{\partial}{\partial t}[\tilde{\zeta}(\omega,t)]\frac{\partial^2}{\partial x^2}[e^{i\lambda x-\omega(\lambda)t}]=e^{\omega(\lambda)t}\zeta(t)+\alpha e^{\omega(\lambda)t}\zeta(t)\lambda^2,$$

the integral in (3.10) becomes

$$\int_{\mathcal{K}(i/\sqrt{\alpha})}\frac{e^{i\lambda x}\zeta(t)}{i\gamma_1\lambda+\gamma_0}\lambda d\lambda,$$

which vanishes by Cauchy's theorem.

*Proof of (3.11)* Since

$$e^{i\lambda x-\omega(\lambda)t}\frac{\partial}{\partial t}[\tilde{\hat{f}}(\lambda,\omega,t)]-\alpha\frac{\partial}{\partial t}[\tilde{\hat{f}}(\lambda,\omega,t)]\frac{\partial^2}{\partial x^2}[e^{i\lambda x-\omega(\lambda)t}]=(1+\alpha\lambda^2)e^{i\lambda x}\hat{f}(\lambda,t)$$

the integral in (3.11) becomes

$$\int_{-\infty}^{\infty}e^{i\lambda x}\hat{f}(\lambda,t)d\lambda=2\pi f(x,t).$$

*Proof of (3.12)* Since

$$e^{i\lambda x-\omega(\lambda)t}\frac{\partial}{\partial t}[\tilde{\hat{f}}(-\lambda,\omega,t)]-\alpha\frac{\partial}{\partial t}[\tilde{\hat{f}}(-\lambda,\omega,t)]\frac{\partial^2}{\partial x^2}[e^{i\lambda x-\omega(\lambda)t}]=(1+\alpha\lambda^2)e^{i\lambda x}\hat{f}(-\lambda,t),$$

the integral in (3.12) becomes

$$\int_{\mathcal{K}(i/\sqrt{\alpha})}e^{i\lambda x}\frac{\gamma_0-i\gamma_1\lambda}{\gamma_0+i\gamma_1\lambda}\hat{f}(-\lambda,t)d\lambda,$$

which clearly vanishes.

***Step 3*** It is clear that

$$\lim_{t\to0^+}\int_{\mathcal{K}(i/\sqrt{\alpha})}\frac{e^{i\lambda x-\omega(\lambda)t}}{(1+\alpha\lambda^2)(i\gamma_1\lambda+\gamma_0)}\lambda\tilde{\zeta}(\omega,t)d\lambda=0 \text{ and } \lim_{t\to0^+}\int_{\mathcal{K}(i/\sqrt{\alpha})}\frac{e^{i\lambda x-\omega(\lambda)t}}{1+\alpha\lambda^2}\frac{\gamma_0-i\gamma_1\lambda}{\gamma_0+i\gamma_1\lambda}\tilde{\hat{f}}(-\lambda,\omega,t)d\lambda=0.$$

Also it is easy to see that

$$\lim_{t\to0^+}\int_{\mathcal{K}(i/\sqrt{\alpha})}\hat{u}_0(-\lambda)\frac{\gamma_0-i\gamma_1\lambda}{\gamma_0+i\gamma_1\lambda}e^{i\lambda x-\omega(\lambda)t}d\lambda=0 \text{ and } \lim_{t\to0^+}\int_{-\infty}^{\infty}\frac{e^{i\lambda x-\omega(\lambda)t}}{1+\alpha\lambda^2}\tilde{\hat{f}}(\lambda,\omega,t)d\lambda=0.$$

Now we claim that

$$\lim_{t\to0^+}\int_{-\infty}^{\infty}\hat{u}_0(\lambda)e^{i\lambda x-\omega(\lambda)t}d\lambda=2\pi u_0(x),\ x>0. \tag{3.13}$$

Since

$$\hat{u}_0(\lambda)=\frac{u_0(0)}{i\lambda}+\frac{\hat{(u_0')}(\lambda)}{i\lambda},\ \lambda\neq0,$$

the integral in (3.13) is equal to

$$\int_{-1}^{1}\hat{u}_0(\lambda)e^{i\lambda x-\omega(\lambda)t}d\lambda+\left(\int_{-\infty}^{-1}+\int_{1}^{\infty}\right)\frac{u_0(0)}{i\lambda}e^{i\lambda x-\omega(\lambda)t}d\lambda+\left(\int_{-\infty}^{-1}+\int_{1}^{\infty}\right)\frac{\hat{(u_0')}(\lambda)}{i\lambda}e^{i\lambda x-\omega(\lambda)t}d\lambda$$

$$=\int_{-1}^{1}\hat{u}_0(\lambda)e^{i\lambda x-\omega(\lambda)t}d\lambda+\int_{\mathcal{C}}\frac{u_0(0)}{i\lambda}e^{i\lambda x-\omega(\lambda)t}d\lambda+\left(\int_{-\infty}^{-1}+\int_{1}^{\infty}\right)\frac{\hat{(u_0')}(\lambda)}{i\lambda}e^{i\lambda x-\omega(\lambda)t}d\lambda.$$

Therefore

$$\lim_{t\to0^+}\int_{-\infty}^{\infty}\hat{u}_0(\lambda)e^{i\lambda x-\omega(\lambda)t}d\lambda=\int_{-1}^{1}\hat{u}_0(\lambda)e^{i\lambda x}d\lambda+\int_{\mathcal{C}}\frac{u_0(0)}{i\lambda}e^{i\lambda x}d\lambda+\left(\int_{-\infty}^{-1}+\int_{1}^{\infty}\right)\frac{\hat{(u_0')}(\lambda)}{i\lambda}e^{i\lambda x}d\lambda$$

$$= \int_{-\infty}^{\infty} \hat{u}_0(\lambda) e^{i\lambda x} d\lambda = 2\pi u_0(x).$$

This proves (3.13).

***Step 4*** We claim that

$$\lim_{x\to 0^+}\left[\mathcal{P}\left(\int_{-\infty}^{\infty} \hat{u}_0(\lambda) e^{i\lambda x-\omega(\lambda)t} d\lambda\right) - \mathcal{P}\left(\int_{\mathcal{K}(i/\sqrt{\alpha})} \hat{u}_0(-\lambda)\frac{\gamma_0 - i\gamma_1\lambda}{\gamma_0 + i\gamma_1\lambda} e^{i\lambda x-\omega(\lambda)t} d\lambda\right)\right] = 0. \tag{3.14}$$

Computing derivatives as in Step 1, we find

$$\left(\alpha\frac{\partial}{\partial t}+\beta\right)\left(\int_{-\infty}^{\infty} \hat{u}_0(\lambda) e^{i\lambda x-\omega(\lambda)t} d\lambda\right) = \beta\int_{-\infty}^{\infty} \hat{u}_0(\lambda)\frac{e^{i\lambda x-\omega(\lambda)t}}{1+\alpha\lambda^2} d\lambda,$$

$$\left(\alpha\frac{\partial}{\partial t}+\beta\right)\left(\int_{-\infty}^{\infty} \hat{u}_0(\lambda) e^{i\lambda x-\omega(\lambda)t} d\lambda\right) = \beta\int_{-\infty}^{\infty} \hat{u}_0(\lambda)\frac{e^{i\lambda x-\omega(\lambda)t}}{1+\alpha\lambda^2} d\lambda,$$

$$\mathcal{P}\left(\int_{-\infty}^{\infty} \hat{u}_0(\lambda) e^{i\lambda x-\omega(\lambda)t} d\lambda\right) = \beta\int_{-\infty}^{\infty} (\gamma_0 + i\gamma_1\lambda)\hat{u}_0(\lambda)\frac{e^{i\lambda x-\omega(\lambda)t}}{1+\alpha\lambda^2} d\lambda,$$

$$\mathcal{P}\left(\int_{\mathcal{K}(i/\sqrt{\alpha})} \hat{u}_0(-\lambda)\frac{\gamma_0 - i\gamma_1\lambda}{\gamma_0 + i\gamma_1\lambda} e^{i\lambda x-\omega(\lambda)t} d\lambda\right) = \beta\int_{\mathcal{K}(i/\sqrt{\alpha})} \hat{u}_0(-\lambda)\frac{\gamma_0 - i\gamma_1\lambda}{1+\alpha\lambda^2} e^{i\lambda x-\omega(\lambda)t} d\lambda$$

$$= \beta\int_{-\infty}^{\infty} \hat{u}_0(-\lambda)\frac{\gamma_0 - i\gamma_1\lambda}{1+\alpha\lambda^2} e^{i\lambda x-\omega(\lambda)t} d\lambda.$$

Thus the limit in (3.14) is equal to

$$\beta\int_{-\infty}^{\infty} (\gamma_0 + i\gamma_1\lambda)\hat{u}_0(\lambda)\frac{e^{-\omega(\lambda)t}}{1+\alpha\lambda^2} d\lambda - \beta\int_{-\infty}^{\infty} \hat{u}_0(-\lambda)\frac{\gamma_0 - i\gamma_1\lambda}{1+\alpha\lambda^2} e^{-\omega(\lambda)t} d\lambda, \tag{3.15}$$

Changing the variable from $\lambda$ to $-\lambda$, in the second integral, we see that (3.15) is indeed zero.

***Step 5*** We claim that

$$\lim_{x\to 0^+}\left[\mathcal{P}\left(\int_{-\infty}^{\infty} \frac{e^{i\lambda x-\omega(\lambda)t}}{1+\alpha\lambda^2}\tilde{\hat{f}}(\lambda,\omega,t) d\lambda\right) - \mathcal{P}\left(\int_{\mathcal{K}(i/\sqrt{\alpha})} \frac{e^{i\lambda x-\omega(\lambda)t}}{1+\alpha\lambda^2}\frac{\gamma_0 - i\gamma_1\lambda}{\gamma_0 + i\gamma_1\lambda}\tilde{\hat{f}}(-\lambda,\omega,t) d\lambda\right)\right] = 0. \tag{3.16}$$

The proof of this is similar to the proof of (3.14), but easier in that the integrands in the integrals of (3.16) are $O(1/\lambda^3)$, and the differentiations of $\mathcal{P}$ can be carried out immediately.

***Step 6*** We claim that

$$\lim_{x\to 0^+}\mathcal{P}\left[\int_{\mathcal{K}(i/\sqrt{\alpha})} \frac{e^{i\lambda x-\omega(\lambda)t}}{(1+\alpha\lambda^2)(i\gamma_1\lambda+\gamma_0)}\lambda\tilde{\zeta}(\omega,t) d\lambda\right] = 2\pi i\zeta(t). \tag{3.17}$$

We compute:

$$\mathcal{P}\left[\int_{\mathcal{K}(i/\sqrt{\alpha})} \cdots d\lambda\right] = \int_{\mathcal{K}(i/\sqrt{\alpha})} \frac{\mathcal{P}[e^{i\lambda x-\omega(\lambda)t}]}{(1+\alpha\lambda^2)(i\gamma_1\lambda+\gamma_0)}\lambda\tilde{\zeta}(\omega,t) d\lambda \tag{3.18}$$

$$+\alpha\gamma_0\int_{\mathcal{K}(i/\sqrt{\alpha})} \frac{e^{i\lambda x-\omega(\lambda)t}}{(1+\alpha\lambda^2)(i\gamma_1\lambda+\gamma_0)}\lambda\frac{\partial}{\partial t}[\tilde{\zeta}(\omega,t)] d\lambda + \alpha\gamma_1\int_{\mathcal{K}(i/\sqrt{\alpha})} \frac{\partial[e^{i\lambda x-\omega(\lambda)t}]/\partial x}{(1+\alpha\lambda^2)(i\gamma_1\lambda+\gamma_0)}\lambda\frac{\partial}{\partial t}[\tilde{\zeta}(\omega,t)] d\lambda$$

$$= \beta \int_{\mathcal{K}(i/\sqrt{\alpha})} \frac{e^{i\lambda x-\omega(\lambda)t}}{(1+\alpha\lambda^2)^2}\lambda\tilde{\zeta}(\omega,t)d\lambda + \alpha h(t)\int_{\mathcal{K}(i/\sqrt{\alpha})}\frac{e^{i\lambda x}}{1+\alpha\lambda^2}\lambda d\lambda := \beta\mathbb{I}_1(x,t)+\alpha\zeta(t)\mathbb{I}_2(x).$$

Now

$$\lim_{x\to 0^+}[\beta\mathbb{I}_1(x,t)] = \beta\int_{\mathcal{K}(i/\sqrt{\alpha})}\frac{e^{-\omega(\lambda)t}}{(1+\alpha\lambda^2)^2}\lambda\tilde{\zeta}(\omega,t)d\lambda = 0. \tag{3.19}$$

In order to prove the second equation in (3.19), we consider the contour $\mathcal{L}(-1/\alpha) := \{\mu=\lambda^2 : \lambda\in\mathcal{K}(i/\sqrt{\alpha})\}$ and we write the integral in the following way:

$$\int_{\mathcal{K}(i/\sqrt{\alpha})}\frac{e^{-\omega(\lambda)t}}{(1+\alpha\lambda^2)^2}\lambda\left[\int_{\tau=0}^{t}e^{\omega(\lambda)\tau}\zeta(\tau)d\tau\right]d\lambda \tag{3.20}$$

$$= \int_{\lambda\in\mathcal{K}(i/\sqrt{\alpha})}\frac{1}{(1+\alpha\lambda^2)^2}\exp\left(-\frac{\beta\lambda^2 t}{1+\alpha\lambda^2}\right)\left[\int_{\tau=0}^{t}\exp\left(\frac{\beta\lambda^2\tau}{1+\alpha\lambda^2}\right)\zeta(\tau)d\tau\right]\lambda d\lambda$$

$$= \frac{1}{2}\int_{\mu\in\mathcal{L}(-1/\alpha)}\frac{1}{(1+\alpha\mu)^2}\exp\left(-\frac{\beta\mu t}{1+\alpha\mu}\right)\left[\int_{\tau=0}^{t}\exp\left(\frac{\beta\mu\tau}{1+\alpha\mu}\right)\zeta(\tau)d\tau\right]d\mu,$$

$$= \frac{1}{2}\int_{\mu\in\mathcal{L}(-1/\alpha)}\frac{1}{(1+\alpha\mu)^2}\left[\int_{\tau=0}^{t}\exp\left(\frac{\beta\mu(t-\tau)}{1+\alpha\mu}\right)\zeta(\tau)d\tau\right]d\mu$$

$$= \frac{1}{2}\int_{\tau=0}^{t}\left[\int_{\mu\in\mathcal{L}(-1/\alpha)}\frac{1}{(1+\alpha\mu)^2}\exp\left(\frac{\beta\mu(t-\tau)}{1+\alpha\mu}\right)d\mu\right]\zeta(\tau)d\tau = 0,$$

where the last equation follows from the facts that

$$\frac{1}{(1+\alpha\mu)^2}\exp\left(\frac{\beta\mu(t-\tau)}{1+\alpha\mu}\right) = \frac{1}{(1+\alpha\mu)^2}\exp\left(\frac{\beta(t-\tau)}{\alpha}\frac{\alpha\mu}{1+\alpha\mu}\right)$$

$$= \frac{1}{(1+\alpha\mu)^2}\exp\left[\frac{\beta(t-\tau)}{\alpha}\left(1-\frac{1}{1+\alpha\mu}\right)\right] = \frac{d}{d\mu}\left\{\frac{1}{\beta(t-\tau)}\exp\left[\frac{\beta(t-\tau)}{\alpha}\left(1-\frac{1}{1+\alpha\mu}\right)\right]\right\}, \text{ for } t\neq\tau,$$

and

$$\frac{1}{(1+\alpha\mu)^2}\exp\left(\frac{\beta\mu(t-\tau)}{1+\alpha\mu}\right) = \frac{d}{d\mu}\left[-\frac{1}{\alpha}\frac{1}{1+\alpha\mu}\right], \text{ for } t=\tau.$$

On the other hand

$$\lim_{x\to 0^+}[\alpha h(t)\mathbb{I}_2(x)] = \alpha\zeta(t)\int_{\mathcal{K}(i/\sqrt{\alpha})}\frac{\lambda d\lambda}{1+\alpha\lambda^2} = \pi i\zeta(t). \tag{3.21}$$

Thus, (3.17) follows from (3.18), (3.19) and (3.20).

***Completion of the proof of Theorem 1*** Equation (1.19) follows from the results of Step 3 and equation (1.20) follows from the results of Steps 4, 5 and 6. The uniformity assertions in the limits of (1.19) and (1.20) follow by examining the proofs in Steps 3 - 6.

***Proof of Theorem 2*** The existence of the limits in the 1st and 2nd parts of the theorem and the assertion of their uniformity follow from (3.5) and (3.9).

It follows from the 1st and 2nd parts of the theorem that the function $u(x,t)$ extends to a $C^\infty$ function $U(x,t)$ for $(x,t)\in\overline{\mathbb{R}^+\times\mathbb{R}^+}-\{(0,0)\}$, by setting

$$U(x,t)=\begin{cases} u_0(x), & for\ t=0\ and\ x>0, \\ \lim_{x\to 0^+} u(x,t), & for\ x=0\ and\ t>0. \end{cases}$$

Finally, the 3rd and 4th parts of the theorem follow from the above conclusion.

***Proof of Theorem 3*** This proof is similar to the proof of the corresponding theorem in [17].

## 4. The problem 2

In this section we will show that the interface *Problem 2* [(1.9), (1.10), (1.11), (1.12)] is separable. As we stated in the *Introduction*, this means that if this problem is solvable, we have to assume (1.13), and this implies that *Problem 2* is reduced to two separate IBV problems, one in the right half-line and another one in the left half-line.

***Derivation of the solution for*** $x>0$ Setting

$$\omega_R=\omega_R(\lambda)=\frac{\beta_R\lambda^2}{1+\alpha_R\lambda^2},\ \alpha_R>0,\ \beta_R>0,$$

we derive the general relation

$$(4.1)\quad \hat{u}^R(\lambda,t)=\hat{u}_0^R(\lambda)e^{-\omega_R(\lambda)t}-\frac{i\lambda e^{-\omega_R(\lambda)t}}{1+\alpha_R\lambda^2}[\alpha_R\tilde{g}'_{0,R}(\omega_R,t)+\beta_R\tilde{g}_{0,R}(\omega_R,t)]$$
$$-\frac{e^{-\omega_R(\lambda)t}}{1+\alpha_R\lambda^2}[\alpha_R\tilde{g}'_{1,R}(\omega_R,t)+\beta_R\tilde{g}_{1,R}(\omega_R,t)]+\frac{e^{-\omega_R(\lambda)t}}{1+\alpha_R\lambda^2}\tilde{\hat{f}}_R(\lambda,\omega_R,t),\ \mathrm{Im}\,\lambda\le 0,$$

and this leads to the formula:

$$(4.2)\quad u^R(x,t)=\frac{1}{2\pi}\int_{-\infty}^{\infty}e^{i\lambda x}\hat{u}^R(\lambda,t)d\lambda$$
$$=\frac{1}{2\pi}\int_{-\infty}^{\infty}\hat{u}_0^R(\lambda)e^{i\lambda x-\omega_R(\lambda)t}d\lambda-\frac{1}{2\pi}\int_{\lambda\in\mathcal{K}(i/\sqrt{\alpha_R})}\frac{i\lambda e^{i\lambda x-\omega_R(\lambda)t}}{1+\alpha_R\lambda^2}[\alpha_R\tilde{g}'_{0,R}(\omega_R,t)+\beta_R\tilde{g}_{0,R}(\omega_R,t)]d\lambda$$
$$-\frac{1}{2\pi}\int_{\lambda\in\mathcal{K}(i/\sqrt{\alpha_R})}\frac{e^{i\lambda x-\omega_R(\lambda)t}}{1+\alpha_R\lambda^2}[\alpha_R\tilde{g}'_{1,R}(\omega_R,t)+\beta_R\tilde{g}_{1,R}(\omega_R,t)]d\lambda+\frac{1}{2\pi}\int_{-\infty}^{\infty}\frac{e^{i\lambda x-\omega_R(\lambda)t}}{1+\alpha_R\lambda^2}\tilde{\hat{f}}_R(\lambda,\omega_R,t)d\lambda,\ x>0,\ t>0.$$

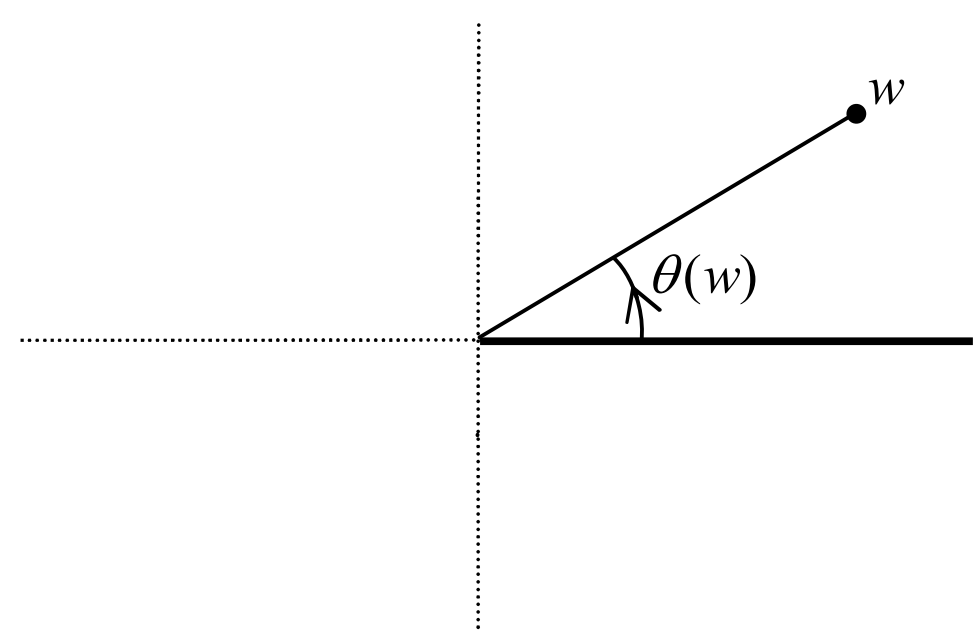

**Fig.3** *The choice of the angle* $\theta(w)$ *in the definition of the square root function* $\sqrt{w}$ : $0\le\theta(w)<2\pi$ .

Now, if

$$\omega_R=\omega_R(\lambda)=\frac{\beta_R\lambda^2}{1+\alpha_R\lambda^2}=\frac{\mu^2}{1+\mu^2}$$
$$\Rightarrow\ \lambda=\varphi_R(\mu):=\sqrt{\frac{\mu^2}{\beta_R+(\beta_R-\alpha_R)\mu^2}}\ \text{ and }\ \mu=\psi_R(\lambda):=\sqrt{\frac{\beta_R\lambda^2}{1+(\alpha_R-\beta_R)\lambda^2}},$$

where the square root is defined by $\sqrt{w} := \sqrt{|w|}e^{i\theta(w)/2}$, $0 \le \theta(w) < 2\pi$. (See fig 3.) This branch of the square root function is analytic for $w \in \mathbb{C} - [0, +\infty)$. Moreover, $\varphi_R(i) = i/\sqrt{\alpha_R}$, $\psi_R(i/\sqrt{\alpha_R}) = i$, and $\varphi_R$ is a biholomorphic mapping from an open neighborhood of $i$ to an open neighborhood of $i/\sqrt{\alpha}$ with inverse $\psi_R$. (See fig 4.)

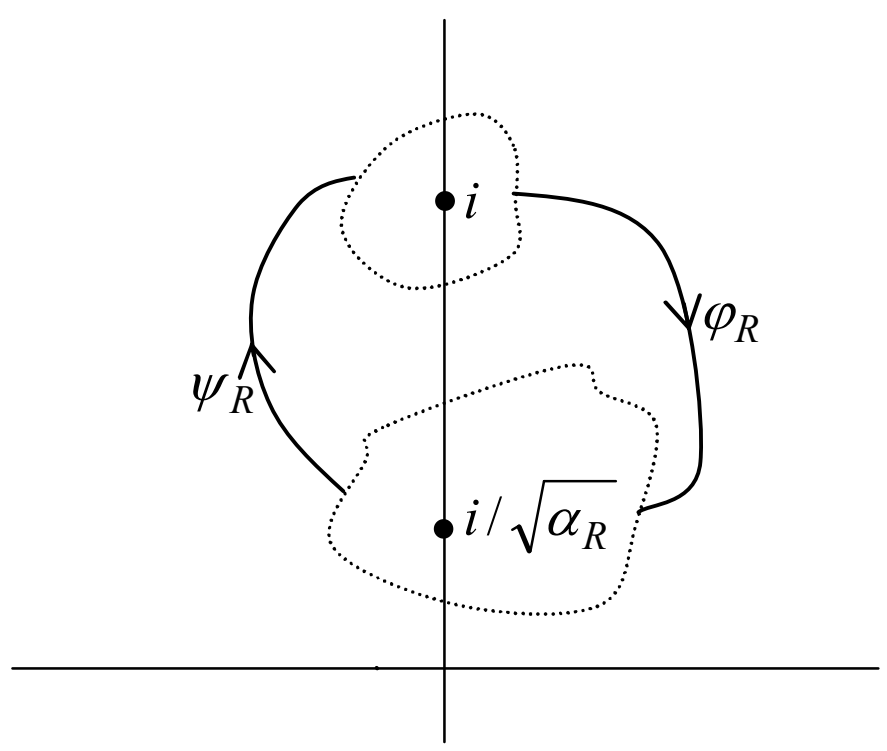

**Fig. 4** $(\varphi_R)'(i) = -\beta_R(\alpha_R)^{-3/2}$ and $(\psi_R)'(i/\sqrt{\alpha_R}) = -(\beta_R)^{-1}(\alpha_R)^{3/2}$.

Setting

$$\omega_0 = \omega_0(\mu) := \frac{\mu^2}{1+\mu^2},$$

and changing the variable in the integrals of (4.2) taken on $\mathcal{K}(i/\sqrt{\alpha_R})$, we obtain that

$$\begin{aligned}(4.3)\; u^R(x,t) &= \frac{1}{2\pi}\int_{-\infty}^{\infty} \hat{u}_0^R(\lambda) e^{i\lambda x - \omega_R(\lambda)t} d\lambda \\ &- \frac{1}{2\pi}\int_{\mu \in C(i,\varepsilon)} \frac{i\varphi_R(\mu) e^{i\varphi_R(\mu)x - \omega_0(\mu)t}}{1+\alpha_R[\varphi_R(\mu)]^2}[\alpha_R \tilde{g}'_{0,R}(\omega_0,t) + \beta_R \tilde{g}_{0,R}(\omega_0,t)]\frac{d[\varphi_R(\mu)]}{d\mu}d\mu \\ &- \frac{1}{2\pi}\int_{\mu \in C(i,\varepsilon)} \frac{e^{i\varphi_R(\mu)x - \omega_0(\mu)t}}{1+\alpha_R[\varphi_R(\mu)]^2}[\alpha_R \tilde{g}'_{1,R}(\omega_0,t) + \beta_R \tilde{g}_{1,R}(\omega_0,t)]\frac{d[\varphi_R(\mu)]}{d\mu}d\mu \\ &+ \frac{1}{2\pi}\int_{-\infty}^{\infty} \frac{e^{i\lambda x - \omega_R(\lambda)t}}{1+\alpha_R\lambda^2}\tilde{\hat{f}}_R(\lambda,\omega_R,t)d\lambda,\; x>0,\; t>0.\end{aligned}$$

***Derivation of the solution for*** $x < 0$ Similarly, setting

$$\omega_L = \omega_L(\lambda) = \frac{\beta_L\lambda^2}{1+\alpha_L\lambda^2},\; \alpha_L > 0, \beta_L > 0,$$

we find that

$$\begin{aligned}(4.4)\quad \hat{u}^L(\lambda,t) &= \hat{u}_0^L(\lambda)e^{-\omega_L(\lambda)t} + \frac{i\lambda e^{-\omega_L(\lambda)t}}{1+\alpha_L\lambda^2}[\alpha_L \tilde{g}'_{0,L}(\omega_L,t) + \beta_L \tilde{g}_{0,L}(\omega_L,t)] \\ &+ \frac{e^{-\omega_L(\lambda)t}}{1+\alpha_L\lambda^2}[\alpha_L \tilde{g}'_{1,L}(\omega_L,t) + \beta_L \tilde{g}_{1,L}(\omega_L,t)] + \frac{e^{-\omega_L(\lambda)t}}{1+\alpha_L\lambda^2}\tilde{\hat{f}}_L(\lambda,\omega_L,t),\; \mathrm{Im}\lambda \ge 0,\end{aligned}$$

and

$$\begin{aligned}(4.5)\quad u^L(x,t) &= \frac{1}{2\pi}\int_{-\infty}^{\infty} e^{i\lambda x}\hat{u}^L(\lambda,t)d\lambda \\ &= \frac{1}{2\pi}\int_{-\infty}^{\infty} \hat{u}_0^L(\lambda)e^{i\lambda x - \omega_L(\lambda)t}d\lambda + \frac{1}{2\pi}\int_{\lambda \in \mathcal{K}(-i/\sqrt{\alpha_L})} \frac{i\lambda e^{i\lambda x - \omega_L(\lambda)t}}{1+\alpha_L\lambda^2}[\alpha_L \tilde{g}'_{0,L}(\omega_L,t) + \beta_L \tilde{g}_{0,L}(\omega_L,t)]d\lambda\end{aligned}$$

$$+\frac{1}{2\pi}\int_{\lambda\in\mathcal{K}(-i/\sqrt{\alpha_L})}\frac{e^{i\lambda x-\omega_R(\lambda)t}}{1+\alpha_L\lambda^2}[\alpha_L\widetilde{g}'_{1,L}(\omega_L,t)+\beta_L\widetilde{g}_{1,L}(\omega_L,t)]d\lambda+\frac{1}{2\pi}\int_{-\infty}^{\infty}\frac{e^{i\lambda x-\omega_L(\lambda)t}}{1+\alpha_L\lambda^2}\widetilde{\hat{f}}_L(\lambda,\omega_L,t)d\lambda,\ x<0,\ t>0,$$

where $\mathcal{K}(-1/\sqrt{\alpha})$ is a sufficiently small simple closed curve surrounding the point $-i/\sqrt{\alpha}$, traversed counterclockwise.

Now, if

$$\omega_L=\omega_L(\lambda)=\frac{\beta_L\lambda^2}{1+\alpha_L\lambda^2}=\frac{\mu^2}{1+\mu^2}$$

$$\Rightarrow\ \lambda=\varphi_L(\mu):=-\sqrt{\frac{\mu^2}{\beta_L+(\beta_L-\alpha_L)\mu^2}}\ \text{ and }\ \mu=\psi_L(\lambda):=\sqrt{\frac{\beta_L\lambda^2}{1+(\alpha_L-\beta_L)\lambda^2}},$$

where the square root is defined as before. Also, $\varphi_L(i)=-i/\sqrt{\alpha}$, $\psi_L(-i/\sqrt{\alpha})=i$, and $\varphi_L$ is a biholomorphic mapping from an open neighborhood of $i$ to an open neighborhood of $-i/\sqrt{\alpha}$ with inverse $\psi_L$. (See fig 5.)

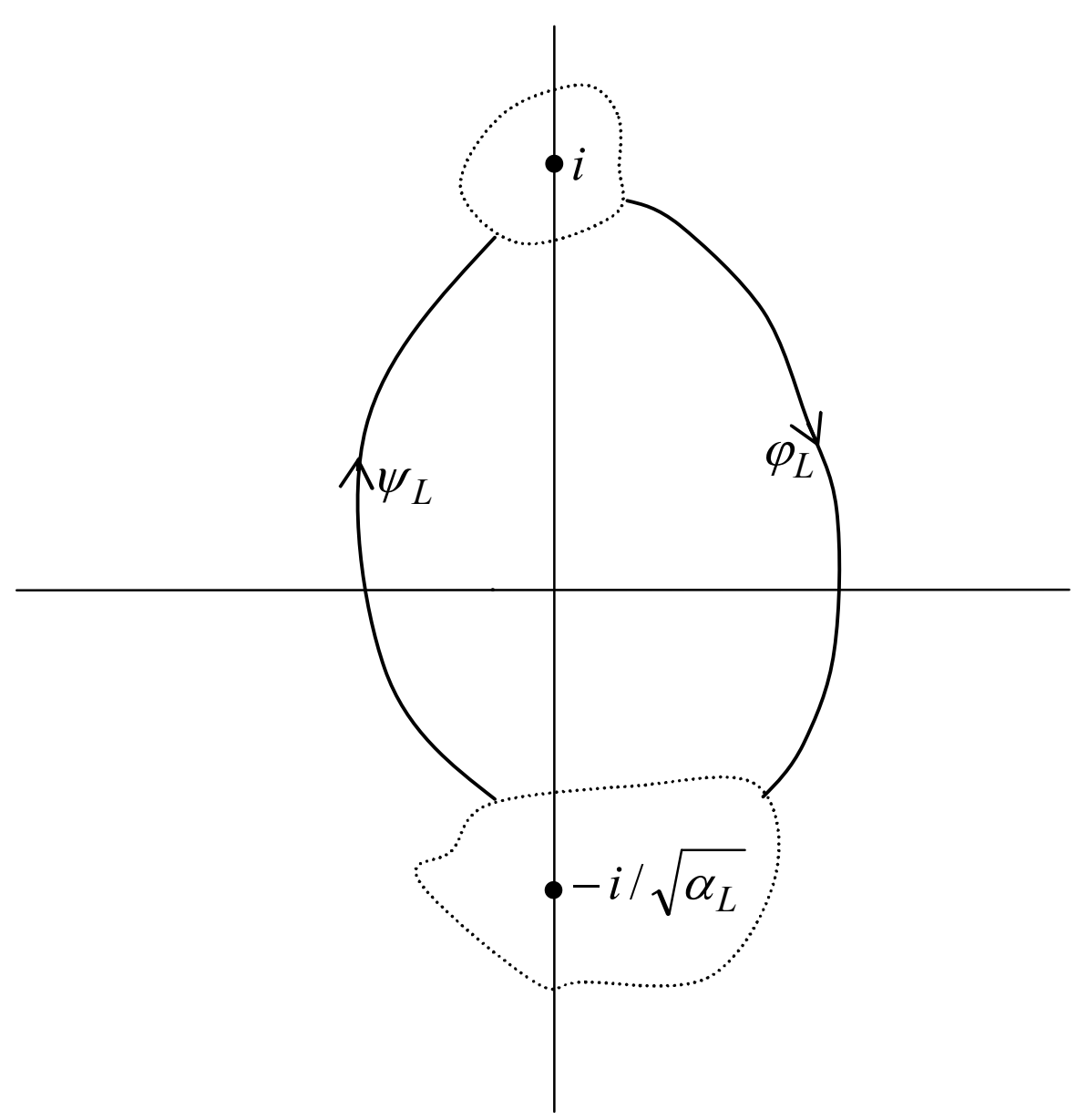

**Fig. 5.** $(\varphi_L)'(i)=\beta_L(\alpha_L)^{-3/2}$ and $(\psi_L)'(-i/\sqrt{\alpha_L})=(\beta_L)^{-1}(\alpha_L)^{3/2}$.

Changing the variable in the integrals of (4.5) taken on $\mathcal{K}(-i/\sqrt{\alpha_L})$, we obtain that

$$\text{(4.6)}\quad u^L(x,t)=\frac{1}{2\pi}\int_{-\infty}^{\infty}\hat{u}_0^L(\lambda)e^{i\lambda x-\omega_L(\lambda)t}d\lambda$$
$$+\frac{1}{2\pi}\int_{\mu\in C(i,\varepsilon)}\frac{i\varphi_L(\mu)e^{i\varphi_L(\mu)x-\omega_0(\mu)t}}{1+\alpha_L[\varphi_L(\mu)]^2}[\alpha_L\widetilde{g}'_{0,L}(\omega_0,t)+\beta_L\widetilde{g}_{0,L}(\omega_0,t)]\frac{d[\varphi_L(\mu)]}{d\mu}d\mu$$
$$+\frac{1}{2\pi}\int_{\mu\in C(i,\varepsilon)}\frac{e^{i\varphi_L(\mu)x-\omega_0(\mu)t}}{1+\alpha_L[\varphi_L(\mu)]^2}[\alpha_L\widetilde{g}'_{1,L}(\omega_0,t)+\beta_L\widetilde{g}_{1,L}(\omega_0,t)]\frac{d[\varphi_L(\mu)]}{d\mu}d\mu$$
$$+\frac{1}{2\pi}\int_{-\infty}^{\infty}\frac{e^{i\lambda x-\omega_L(\lambda)t}}{1+\alpha_L\lambda^2}\widetilde{\hat{f}}_L(\lambda,\omega_L,t)d\lambda,\ x<0,\ t>0.$$

***The attempt to compute the integrals containing unknown quantities*** The unknown quantities we are referring to are $\hat{u}^R(-\varphi_R(\mu),t)$ and $\hat{u}^L(-\varphi_L(\mu),t)$.

Setting $-\lambda$ in (4.1) and (4.4) we obtain

$$(4.7)\quad \hat{u}^R(-\lambda,t)=\hat{u}_0^R(-\lambda)e^{-\omega_R(\lambda)t}+\frac{i\lambda e^{-\omega_R(\lambda)t}}{1+\alpha_R\lambda^2}[\alpha_R\tilde{g}'_{0,R}(\omega_R,t)+\beta_R\tilde{g}_{0,R}(\omega_R,t)]$$
$$-\frac{e^{-\omega_R(\lambda)t}}{1+\alpha_R\lambda^2}[\alpha_R\tilde{g}'_{1,R}(\omega_R,t)+\beta_R\tilde{g}_{1,R}(\omega_R,t)]+\frac{e^{-\omega_R(\lambda)t}}{1+\alpha_R\lambda^2}\tilde{\hat{f}}_R(-\lambda,\omega_R,t),\ \operatorname{Im}\lambda\geq 0,$$

and

$$(4.8)\quad \hat{u}^L(-\lambda,t)=\hat{u}_0^L(-\lambda)e^{-\omega_L(\lambda)t}-\frac{i\lambda e^{-\omega_L(\lambda)t}}{1+\alpha_L\lambda^2}[\alpha_L\tilde{g}'_{0,L}(\omega_L,t)+\beta_L\tilde{g}_{0,L}(\omega_L,t)]$$
$$+\frac{e^{-\omega_L(\lambda)t}}{1+\alpha_L\lambda^2}[\alpha_L\tilde{g}'_{1,L}(\omega_L,t)+\beta_L\tilde{g}_{1,L}(\omega_L,t)]+\frac{e^{-\omega_L(\lambda)t}}{1+\alpha_L\lambda^2}\tilde{\hat{f}}_L(-\lambda,\omega_L,t),\ \operatorname{Im}\lambda\leq 0.$$

Then (4.7) give

$$(4.9)\quad \hat{u}^R(-\varphi_R(\mu),t)=\hat{u}_0^R(-\varphi_R(\mu))e^{-\omega_0(\mu)t}+\frac{i\varphi_R(\mu)e^{-\omega_0(\mu)t}}{1+\alpha_R[\varphi_R(\mu)]^2}[\alpha_R\tilde{g}'_{0,R}(\omega_0,t)+\beta_R\tilde{g}_{0,R}(\omega_0,t)]$$
$$-\frac{e^{-\omega_0(\mu)t}}{1+\alpha_R[\varphi_R(\mu)]^2}[\alpha_R\tilde{g}'_{1,R}(\omega_0,t)+\beta_R\tilde{g}_{1,R}(\omega_0,t)]+\frac{e^{-\omega_0(\mu)t}}{1+\alpha_R[\varphi_R(\mu)]^2}\tilde{\hat{f}}_R(-\varphi_R(\mu),\omega_0,t),\ \operatorname{Im}\varphi_R(\mu)\geq 0.$$

In particular (4.9) holds for $\mu\in C(i,\varepsilon)$.

Also (4.8) give

$$(4.10)\quad \hat{u}^L(-\varphi_L(\mu),t)=\hat{u}_0^L(-\varphi_L(\mu))e^{-\omega_0(\mu)t}-\frac{i\varphi_L(\mu)e^{-\omega_0(\mu)t}}{1+\alpha_L[\varphi_L(\mu)]^2}[\alpha_L\tilde{g}'_{0,L}(\omega_0,t)+\beta_L\tilde{g}_{0,L}(\omega_0,t)]$$
$$+\frac{e^{-\omega_0(\mu)t}}{1+\alpha_L[\varphi_L(\mu)]^2}[\alpha_L\tilde{g}'_{1,L}(\omega_0,t)+\beta_L\tilde{g}_{1,L}(\omega_0,t)]+\frac{e^{-\omega_0(\mu)t}}{1+\alpha_L[\varphi_L(\mu)]^2}\tilde{\hat{f}}_L(-\varphi_L(\mu),\omega_0,t),\ \operatorname{Im}\varphi_L(\mu)\leq 0.$$

In particular (4.10) holds for $\mu\in C(i,\varepsilon)$.

Setting

$$G^*_{0,R}(\omega_0,t):=\alpha_R\tilde{g}'_{0,R}(\omega_0,t)+\beta_R\tilde{g}_{0,R}(\omega_0,t),\ G^*_{1,R}(\omega_0,t):=\alpha_R\tilde{g}'_{1,R}(\omega_0,t)+\beta_R\tilde{g}_{1,R}(\omega_0,t),$$
$$G^*_{0,L}(\omega_0,t):=\alpha_L\tilde{g}'_{0,L}(\omega_0,t)+\beta_L\tilde{g}_{0,L}(\omega_0,t),\ G^*_{1,L}(\omega_0,t):=\alpha_L\tilde{g}'_{1,L}(\omega_0,t)+\beta_L\tilde{g}_{1,L}(\omega_0,t),$$

we write (4.9) and (4.10) as follows: For $\mu\in C(i,\varepsilon)$,

$$i\varphi_R(\mu)G^*_{0,R}(\omega_0,t)-G^*_{1,R}(\omega_0,t)=\{1+\alpha_R[\varphi_R(\mu)]^2\}\{e^{\omega_0(\mu)t}\hat{u}^R(-\varphi_R(\mu),t)-\hat{u}_0^R(-\varphi_R(\mu))\}-\tilde{\hat{f}}_R(-\varphi_R(\mu),\omega_0,t),$$
$$i\varphi_L(\mu)G^*_{0,L}(\omega_0,t)-G^*_{1,L}(\omega_0,t)=-\{1+\alpha_L[\varphi_L(\mu)]^2\}\{e^{\omega_0(\mu)t}\hat{u}^L(-\varphi_L(\mu),t)+\hat{u}_0^L(-\varphi_L(\mu))\}+\tilde{\hat{f}}_L(-\varphi_L(\mu),\omega_0,t).$$

Also,

$$\gamma_{11}[\alpha_R\tilde{g}'_{0,R}(\omega_0,t)+\beta_R\tilde{g}_{0,R}(\omega_0,t)]+\gamma_{12}[\alpha_R\tilde{g}'_{1,R}(\omega_0,t)+\beta_R\tilde{g}_{1,R}(\omega_0,t)]$$
$$+\gamma_{13}[\alpha_L\tilde{g}'_{0,L}(\omega_0,t)+\beta_L\tilde{g}_{0,L}(\omega_0,t)]+\gamma_{14}[\alpha_L\tilde{g}'_{1,L}(\omega_0,t)+\beta_L\tilde{g}_{1,L}(\omega_0,t)]=\tilde{I}_1(\omega_0,t),$$

$$\gamma_{21}[\alpha_R\tilde{g}'_{0,R}(\omega_0,t)+\beta_R\tilde{g}_{0,R}(\omega_0,t)]+\gamma_{22}[\alpha_R\tilde{g}'_{1,R}(\omega_0,t)+\beta_R\tilde{g}_{1,R}(\omega_0,t)]$$
$$+\gamma_{23}[\alpha_L\tilde{g}'_{0,L}(\omega_0,t)+\beta_L\tilde{g}_{0,L}(\omega_0,t)]+\gamma_{24}[\alpha_L\tilde{g}'_{1,L}(\omega_0,t)+\beta_L\tilde{g}_{1,L}(\omega_0,t)]=\tilde{I}_2(\omega_0,t),$$

and with the notation we introduced, the above equations can be written as follows:

$$(4.11)\quad \gamma_{11}G^*_{0,R}(\omega_0,t)+\gamma_{12}G^*_{1,R}(\omega_0,t)+\gamma_{13}G^*_{0,L}(\omega_0,t)+\gamma_{14}G^*_{1,L}(\omega_0,t)=\tilde{I}_1(\omega_0,t),$$

(4.12) $\gamma_{21}G^{*}_{0,R}(\omega_0,t)+\gamma_{22}G^{*}_{1,R}(\omega_0,t)+\gamma_{23}G^{*}_{0,L}(\omega_0,t)+\gamma_{24}G^{*}_{1,L}(\omega_0,t)=\widetilde{I}_2(\omega_0,t)$,

and these hold for every $\mu\in\mathbb{C}$.

Now, we procced to solve the system of equation (4.9), (4.10), (4.11) and (4.12), with $\mu\in C(i,\varepsilon)$, in the unknowns

$$G^{*}_{0,R}(\omega_0,t),\, G^{*}_{1,R}(\omega_0,t),\, G^{*}_{0,L}(\omega_0,t),\, G^{*}_{1,L}(\omega_0,t).$$

By Cramer's rule, we have:

$$G^{*}_{0,R}(\mu,t)=\frac{\mathfrak{D}^{R,*}_0(\mu,t)}{\mathfrak{D}(\mu)},\; G^{*}_{1,R}(\mu,t)=\frac{\mathfrak{D}^{R,*}_1(\mu,t)}{\mathfrak{D}(\mu)},\; G^{*}_{0,L}(\mu,t)=\frac{\mathfrak{D}^{L,*}_0(\mu,t)}{\mathfrak{D}(\mu)},\; G^{*}_{1,L}(\mu,t)=\frac{\mathfrak{D}^{L,*}_1(\mu,t)}{\mathfrak{D}(\mu)},\; \mu\in C(i,\varepsilon),$$

where we used the following notation:

$$\mathfrak{D}(\mu):=\det\begin{bmatrix} i\varphi_R(\mu) & -1 & 0 & 0 \\ 0 & 0 & i\varphi_L(\mu) & -1 \\ \gamma_{11} & \gamma_{12} & \gamma_{13} & \gamma_{14} \\ \gamma_{21} & \gamma_{22} & \gamma_{23} & \gamma_{24} \end{bmatrix},$$

$$\mathfrak{D}^{R,*}_0(\mu,t):=\det\begin{bmatrix} J^{*}_R(\mu,t)+J_R(\mu) & -1 & 0 & 0 \\ J^{*}_L(\mu,t)+J_L(\mu) & 0 & i\varphi_L(\mu) & -1 \\ \widetilde{I}_1(\omega_0,t) & \gamma_{12} & \gamma_{13} & \gamma_{14} \\ \widetilde{I}_2(\omega_0,t) & \gamma_{22} & \gamma_{23} & \gamma_{24} \end{bmatrix},$$

$$\mathfrak{D}^{R,*}_1(\mu,t):=\det\begin{bmatrix} i\varphi_R(\mu) & J^{*}_R(\mu,t)+J_R(\mu) & 0 & 0 \\ 0 & J^{*}_L(\mu,t)+J_L(\mu) & i\varphi_L(\mu) & -1 \\ \gamma_{11} & \widetilde{I}_1(\omega_0,t) & \gamma_{13} & \gamma_{14} \\ \gamma_{21} & \widetilde{I}_2(\omega_0,t) & \gamma_{23} & \gamma_{24} \end{bmatrix},$$

$$\mathfrak{D}^{L,*}_0(\mu,t):=\det\begin{bmatrix} i\varphi_R(\mu) & -1 & J^{*}_R(\mu,t)+J_R(\mu) & 0 \\ 0 & 0 & J^{*}_L(\mu,t)+J_L(\mu) & -1 \\ \gamma_{11} & \gamma_{12} & \widetilde{I}_1(\omega_0,t) & \gamma_{14} \\ \gamma_{21} & \gamma_{22} & \widetilde{I}_2(\omega_0,t) & \gamma_{24} \end{bmatrix},$$

$$\mathfrak{D}^{L,*}_1(\mu,t):=\det\begin{bmatrix} i\varphi_R(\mu) & -1 & 0 & J^{*}_R(\mu,t)+J_R(\mu) \\ 0 & 0 & i\varphi_L(\mu) & J^{*}_L(\mu,t)+J_L(\mu) \\ \gamma_{11} & \gamma_{12} & \gamma_{13} & \widetilde{I}_1(\omega_0,t) \\ \gamma_{21} & \gamma_{22} & \gamma_{23} & \widetilde{I}_2(\omega_0,t) \end{bmatrix},$$

$$J^{*}_R(\mu,t):=\{1+\alpha_R[\varphi_R(\mu)]^2\}e^{\omega_0(\mu)t}\hat{u}^R(-\varphi_R(\mu),t),\;\; J_R(\mu,t):=-\{1+\alpha_R[\varphi_R(\mu)]^2\}\hat{u}^R_0(-\varphi_R(\mu))-\widetilde{\widetilde{f}}_R(-\varphi_R(\mu),\omega_0,t),$$

$$J^{*}_L(\mu,t):=-\{1+\alpha_L[\varphi_L(\mu)]^2\}e^{\omega_0(\mu)t}\hat{u}^L(-\varphi_L(\mu),t),\;\; J_L(\mu,t):=\{1+\alpha_L[\varphi_L(\mu)]^2\}\hat{u}^L_0(-\varphi_L(\mu))-\widetilde{\widetilde{f}}_L(-\varphi_L(\mu),\omega_0,t).$$

Let us define also

$$G_{0,R}(\mu,t):=\frac{\mathfrak{D}^{R}_0(\mu,t)}{\mathfrak{D}(\mu)},\; G_{1,R}(\mu,t):=\frac{\mathfrak{D}^{R}_1(\mu,t)}{\mathfrak{D}(\mu)},\; G_{0,L}(\mu,t):=\frac{\mathfrak{D}^{L}_0(\mu,t)}{\mathfrak{D}(\mu)},\; G_{1,L}(\mu,t):=\frac{\mathfrak{D}^{L}_1(\mu,t)}{\mathfrak{D}(\mu)},$$

where $\mathfrak{D}_0^R(\mu,t)$, $\mathfrak{D}_1^R(\mu,t)$, $\mathfrak{D}_0^L(\mu,t)$, $\mathfrak{D}_1^L(\mu,t)$ are the deterninants which result from the corresponding deterninants $\mathfrak{D}_0^{R,*}(\mu,t)$, $\mathfrak{D}_1^{R,*}(\mu,t)$, $\mathfrak{D}_0^{L,*}(\mu,t)$, $\mathfrak{D}_1^{L,*}(\mu,t)$, by eliminating the terms $J_R^*(\mu,t)$ and $J_L^*(\mu,t)$, i.e., by replacing them by zero, i.e.,

$$\mathfrak{D}_0^R(\mu,t):=\det\begin{bmatrix} J_R(\mu) & -1 & 0 & 0 \\ J_L(\mu) & 0 & i\varphi_L(\mu) & -1 \\ \widetilde{I}_1(\omega_0,t) & \gamma_{12} & \gamma_{13} & \gamma_{14} \\ \widetilde{I}_2(\omega_0,t) & \gamma_{22} & \gamma_{23} & \gamma_{24} \end{bmatrix},\ \mathfrak{D}_1^R(\mu,t):=\det\begin{bmatrix} i\varphi_R(\mu) & J_R(\mu) & 0 & 0 \\ 0 & J_L(\mu) & i\varphi_L(\mu) & -1 \\ \gamma_{11} & \widetilde{I}_1(\omega_0,t) & \gamma_{13} & \gamma_{14} \\ \gamma_{21} & \widetilde{I}_2(\omega_0,t) & \gamma_{23} & \gamma_{24} \end{bmatrix},$$

$$\mathfrak{D}_0^L(\mu,t):=\det\begin{bmatrix} i\varphi_R(\mu) & -1 & J_R(\mu) & 0 \\ 0 & 0 & J_L(\mu) & -1 \\ \gamma_{11} & \gamma_{12} & \widetilde{I}_1(\omega_0,t) & \gamma_{14} \\ \gamma_{21} & \gamma_{22} & \widetilde{I}_2(\omega_0,t) & \gamma_{24} \end{bmatrix},\mathfrak{D}_1^L(\mu,t):=\det\begin{bmatrix} i\varphi_R(\mu) & -1 & 0 & J_R(\mu) \\ 0 & 0 & i\varphi_L(\mu) & J_L(\mu) \\ \gamma_{11} & \gamma_{12} & \gamma_{13} & \widetilde{I}_1(\omega_0,t) \\ \gamma_{21} & \gamma_{22} & \gamma_{23} & \widetilde{I}_2(\omega_0,t) \end{bmatrix}.$$

Thus,

$$\mathfrak{D}_0^{R,*}(\mu,t)-\mathfrak{D}_0^R(\mu,t)=\mathfrak{U}_0^R(\mu,t),\ \ \mathfrak{D}_1^{R,*}(\mu,t)-\mathfrak{D}_1^R(\mu,t)=\mathfrak{U}_1^R(\mu,t),$$

$$\mathfrak{D}_0^{L,*}(\mu,t)-\mathfrak{D}_0^L(\mu,t)=\mathfrak{U}_0^L(\mu,t),\ \ \mathfrak{D}_1^{L,*}(\mu,t)-\mathfrak{D}_1^L(\mu,t)=\mathfrak{U}_1^L(\mu,t),$$

where

$$\mathfrak{U}_0^R(\mu,t):=\det\begin{bmatrix} J_R^*(\mu,t) & -1 & 0 & 0 \\ J_L^*(\mu,t) & 0 & i\varphi_L(\mu) & -1 \\ 0 & \gamma_{12} & \gamma_{13} & \gamma_{14} \\ 0 & \gamma_{22} & \gamma_{23} & \gamma_{24} \end{bmatrix},$$

$$\mathfrak{U}_1^R(\mu,t):=\det\begin{bmatrix} i\varphi_R(\mu) & J_R^*(\mu,t) & 0 & 0 \\ 0 & J_L^*(\mu,t) & i\varphi_L(\mu) & -1 \\ \gamma_{11} & 0 & \gamma_{13} & \gamma_{14} \\ \gamma_{21} & 0 & \gamma_{23} & \gamma_{24} \end{bmatrix},$$

$$\mathfrak{U}_0^L(\mu,t):=\det\begin{bmatrix} i\varphi_R(\mu) & -1 & J_R^*(\mu,t) & 0 \\ 0 & 0 & J_L^*(\mu,t) & -1 \\ \gamma_{11} & \gamma_{12} & 0 & \gamma_{14} \\ \gamma_{21} & \gamma_{22} & 0 & \gamma_{24} \end{bmatrix},$$

$$\mathfrak{U}_1^L(\mu,t):=\det\begin{bmatrix} i\varphi_R(\mu) & -1 & 0 & J_R^*(\mu,t) \\ 0 & 0 & i\varphi_L(\mu) & J_L^*(\mu,t) \\ \gamma_{11} & \gamma_{12} & \gamma_{13} & 0 \\ \gamma_{21} & \gamma_{22} & \gamma_{23} & 0 \end{bmatrix}.$$

Now according to the unified transform method the integrals

$$(4.13)\qquad \int\limits_{\mu\in C(i,\varepsilon)} \frac{i\varphi_R(\mu)e^{i\varphi_R(\mu)x-\omega_0(\mu)t}}{1+\alpha_R[\varphi_R(\mu)]^2}\mathfrak{U}_0^R(\mu,t)\frac{d[\varphi_R(\mu)]}{d\mu}d\mu\,,\quad \int\limits_{\mu\in C(i,\varepsilon)} \frac{e^{i\varphi_R(\mu)x-\omega_0(\mu)t}}{1+\alpha_R[\varphi_R(\mu)]^2}\mathfrak{U}_1^R(\mu,t)\frac{d[\varphi_R(\mu)]}{d\mu}d\mu$$

and

(4.14) $$\int_{\mu\in C(i,\varepsilon)} \frac{i\varphi_L(\mu)e^{i\varphi_L(\mu)x-\omega_0(\mu)t}}{1+\alpha_L[\varphi_L(\mu)]^2}\mathfrak{U}_0^L(\mu,t)\frac{d[\varphi_L(\mu)]}{d\mu}d\mu\,,\quad \int_{\mu\in C(i,\varepsilon)} \frac{e^{i\varphi_L(\mu)x-\omega_0(\mu)t}}{1+\alpha_L[\varphi_L(\mu)]^2}\mathfrak{U}_1^L(\mu,t)\frac{d[\varphi_L(\mu)]}{d\mu}d\mu$$

must vanish, since these integrals contain the unknown quantities $J_R^*(\mu,t)$, $J_L^*(\mu,t)$. The parts of the integrals (4.13) which contain the term $J_R^*(\mu,t)$ do vanish, but this is not true with the parts of these integrals, namely (4.13), which contain the term $J_L^*(\mu,t)$ (*unless* $\alpha_R=\alpha_L$). Thus, in order for the UTM to be effective we have to assume that the coefficients of $J_L^*(\mu,t)$ in the terms $\mathfrak{U}_k^R(\mu,t)$, $k=0,1$, are zero. This leads to the condition $\Gamma_{34}=0$.

Similarly, in order to ensure the vanishing of the parts of the integrals (4.14) which contain the term $J_R^*(\mu,t)$, we have to assume the condition $\Gamma_{12}=0$.

## 5. A general interface problem

***Problem*** *3 Solve the initial value problems (1.9) and (1.10), subject to the following interface conditions:*

*For* $t>0$*,*

(5.1) $$\begin{cases}\gamma_{11}u^R(0,t)+\gamma_{12}u_x^R(0,t)+\gamma_{13}u_t^R(0,t)+\gamma_{14}u_{xt}^R(0,t)+\gamma_{15}u^L(0,t)+\gamma_{16}u_x^L(0,t)+\gamma_{17}u_t^L(0,t)+\gamma_{18}u_{xt}^L(0,t)=I_1(t),\\ \gamma_{21}u^R(0,t)+\gamma_{22}u_x^R(0,t)+\gamma_{23}u_t^R(0,t)+\gamma_{24}u_{xt}^R(0,t)+\gamma_{25}u^L(0,t)+\gamma_{26}u_x^L(0,t)+\gamma_{27}u_t^L(0,t)+\gamma_{28}u_{xt}^L(0,t)=I_2(t),\\ \gamma_{31}u^R(0,t)+\gamma_{32}u_x^R(0,t)+\gamma_{33}u_t^R(0,t)+\gamma_{34}u_{xt}^R(0,t)+\gamma_{35}u^L(0,t)+\gamma_{36}u_x^L(0,t)+\gamma_{37}u_t^L(0,t)+\gamma_{38}u_{xt}^L(0,t)=I_3(t),\\ \gamma_{41}u^R(0,t)+\gamma_{42}u_x^R(0,t)+\gamma_{43}u_t^R(0,t)+\gamma_{44}u_{xt}^R(0,t)+\gamma_{45}u^L(0,t)+\gamma_{46}u_x^L(0,t)+\gamma_{47}u_t^L(0,t)+\gamma_{48}u_{xt}^L(0,t)=I_4(t),\\ \gamma_{51}u^R(0,t)+\gamma_{52}u_x^R(0,t)+\gamma_{53}u_t^R(0,t)+\gamma_{54}u_{xt}^R(0,t)+\gamma_{55}u^L(0,t)+\gamma_{56}u_x^L(0,t)+\gamma_{57}u_t^L(0,t)+\gamma_{58}u_{xt}^L(0,t)=I_5(t),\\ \gamma_{61}u^R(0,t)+\gamma_{62}u_x^R(0,t)+\gamma_{63}u_t^R(0,t)+\gamma_{64}u_{xt}^R(0,t)+\gamma_{65}u^L(0,t)+\gamma_{66}u_x^L(0,t)+\gamma_{67}u_t^L(0,t)+\gamma_{68}u_{xt}^L(0,t)=I_6(t).\end{cases}$$

***Derivation of the solution*** Working as in section 4, we derive the general relations (4.1) and (4.4), the first representation formulas of the solutions (4.2) and (4.5). Continuing as in section 4 and using the maps $\varphi_R(\mu)$ and $\varphi_L(\mu)$, we obtain the representation formulas (4.3) and (4.6).

Now – and this step is different from what we did in section 4 – we write the equations:

(5.2) $$\begin{aligned}&\gamma_{11}\tilde{g}_{0,R}+\gamma_{12}\tilde{g}_{1,R}+\gamma_{13}\tilde{g}'_{0,R}+\gamma_{14}\tilde{g}'_{1,R}+\gamma_{15}\tilde{g}_{0,L}+\gamma_{16}\tilde{g}_{1,L}+\gamma_{17}\tilde{g}'_{0,L}+\gamma_{18}\tilde{g}'_{1,L}=\tilde{I}_1(\omega_0,t)\\ &\gamma_{21}\tilde{g}_{0,R}+\gamma_{22}\tilde{g}_{1,R}+\gamma_{23}\tilde{g}'_{0,R}+\gamma_{24}\tilde{g}'_{1,R}+\gamma_{25}\tilde{g}_{0,L}+\gamma_{26}\tilde{g}_{1,L}+\gamma_{27}\tilde{g}'_{0,L}+\gamma_{28}\tilde{g}'_{1,L}=\tilde{I}_2(\omega_0,t)\\ &\gamma_{31}\tilde{g}_{0,R}+\gamma_{32}\tilde{g}_{1,R}+\gamma_{33}\tilde{g}'_{0,R}+\gamma_{34}\tilde{g}'_{1,R}+\gamma_{35}\tilde{g}_{0,L}+\gamma_{36}\tilde{g}_{1,L}+\gamma_{37}\tilde{g}'_{0,L}+\gamma_{38}\tilde{g}'_{1,L}=\tilde{I}_3(\omega_0,t)\\ &\gamma_{41}\tilde{g}_{0,R}+\gamma_{42}\tilde{g}_{1,R}+\gamma_{43}\tilde{g}'_{0,R}+\gamma_{44}\tilde{g}'_{1,R}+\gamma_{45}\tilde{g}_{0,L}+\gamma_{46}\tilde{g}_{1,L}+\gamma_{47}\tilde{g}'_{0,L}+\gamma_{48}\tilde{g}'_{1,L}=\tilde{I}_4(\omega_0,t)\\ &\gamma_{51}\tilde{g}_{0,R}+\gamma_{52}\tilde{g}_{1,R}+\gamma_{53}\tilde{g}'_{0,R}+\gamma_{54}\tilde{g}'_{1,R}+\gamma_{55}\tilde{g}_{0,L}+\gamma_{56}\tilde{g}_{1,L}+\gamma_{57}\tilde{g}'_{0,L}+\gamma_{58}\tilde{g}'_{1,L}=\tilde{I}_5(\omega_0,t)\\ &\gamma_{61}\tilde{g}_{0,R}+\gamma_{62}\tilde{g}_{1,R}+\gamma_{63}\tilde{g}'_{0,R}+\gamma_{64}\tilde{g}'_{1,R}+\gamma_{65}\tilde{g}_{0,L}+\gamma_{66}\tilde{g}_{1,L}+\gamma_{67}\tilde{g}'_{0,L}+\gamma_{68}\tilde{g}'_{1,L}=\tilde{I}_6(\omega_0,t)\,,\end{aligned}$$

which follow from the interface conditions (5.1).

Next we solve the system of equations (5.2), (4.9) and (4.10), in the unknowns

$\tilde{g}_{0,R}=\tilde{g}_{0,R}(\omega_0(\mu),t)$, $\tilde{g}_{1,R}=\tilde{g}_{1,R}(\omega_0(\mu),t)$, $\tilde{g}'_{0,R}=\tilde{g}'_{0,R}(\omega_0(\mu),t)$, $\tilde{g}'_{1,R}=\tilde{g}'_{1,R}(\omega_0(\mu),t)$,

$\tilde{g}_{0,L}=\tilde{g}_{0,L}(\omega_0(\mu),t)$, $\tilde{g}_{1,L}=\tilde{g}_{1,L}(\omega_0(\mu),t)$, $\tilde{g}'_{0,L}=\tilde{g}'_{0,L}(\omega_0(\mu),t)$, $\tilde{g}'_{1,L}=\tilde{g}'_{1,L}(\omega_0(\mu),t)$,

with $\mu\in C(i,\varepsilon)$, and for fixed $t>0$.

Let us write (4.9) and (4.10) in the form:

(5.3) $$i\varphi_R(\mu)[\alpha_R\tilde{g}'_{0,R}(\omega_0,t)+\beta_R\tilde{g}_{0,R}(\omega_0,t)]-[\alpha_R\tilde{g}'_{1,R}(\omega_0,t)+\beta_R\tilde{g}_{1,R}(\omega_0,t)]=J^*_R(\mu,t)+J_R(\mu,t)\,,\ \mu\in C(i,\varepsilon)\,,$$

(5.4) $$i\varphi_L(\mu)[\alpha_L\tilde{g}'_{0,L}(\omega_0,t)+\beta_L\tilde{g}_{0,L}(\omega_0,t)]-[\alpha_L\tilde{g}'_{1,L}(\omega_0,t)+\beta_L\tilde{g}_{1,L}(\omega_0,t)]=J^*_L(\mu,t)+J_L(\mu,t)\,,\ \mu\in C(i,\varepsilon)\,,$$

where we have set

$$J^*_R(\mu,t):=\{1+\alpha_R[\varphi_R(\mu)]^2\}e^{\omega_0(\mu)t}\hat{u}^R(-\varphi_R(\mu),t)\,,\ J_R(\mu,t):=-\{1+\alpha_R[\varphi_R(\mu)]^2\}\hat{u}_0^R(-\varphi_R(\mu))-\tilde{\hat{f}}_R(-\varphi_R(\mu),\omega_0,t)$$

and

$$J^*_L(\mu,t):=-\{1+\alpha_L[\varphi_L(\mu)]^2\}e^{\omega_0(\mu)t}\hat{u}^L(-\varphi_L(\mu),t)\,,\ J_L(\mu,t):=\{1+\alpha_L[\varphi_L(\mu)]^2\}\hat{u}_0^L(-\varphi_L(\mu))-\tilde{\hat{f}}_L(-\varphi_L(\mu),\omega_0,t)\,.$$

Solving the system of equations (5.2), 5.3) and (5.4), we find

(5.5) $$\tilde{g}_{0,R}=\frac{\det G^*_1(\mu,t)}{\det[\mathfrak{Coef}(\mu)]}\,,\ \tilde{g}_{1,R}=\frac{\det G^*_2(\mu,t)}{\det[\mathfrak{Coef}(\mu)]}\,,\ \tilde{g}'_{0,R}=\frac{\det G^*_3(\mu,t)}{\det[\mathfrak{Coef}(\mu)]}\,,\ \tilde{g}'_{1,R}=\frac{\det G^*_4(\mu,t)}{\det[\mathfrak{Coef}(\mu)]}\,,$$

(5.6) $$\tilde{g}_{0,L}=\frac{\det G^*_5(\mu,t)}{\det[\mathfrak{Coef}(\mu)]}\,,\ \tilde{g}_{1,L}=\frac{\det G^*_6(\mu,t)}{\det[\mathfrak{Coef}(\mu)]}\,,\ \tilde{g}'_{0,L}=\frac{\det G^*_7(\mu,t)}{\det[\mathfrak{Coef}(\mu)]}\,,\ \tilde{g}'_{1,L}=\frac{\det G^*_8(\mu,t)}{\det[\mathfrak{Coef}(\mu)]}\,,$$

where we use the following notation:

$$\mathfrak{Coef}(\mu):=\begin{bmatrix}\beta_R i\varphi_R(\mu) & -\beta_R & \alpha_R i\varphi_R(\mu) & -\alpha_R & 0 & 0 & 0 & 0\\ 0 & 0 & 0 & 0 & \beta_L i\varphi_L(\mu) & -\beta_L & \alpha_L i\varphi_L(\mu) & -\alpha_L\\ \gamma_{11} & \gamma_{12} & \gamma_{13} & \gamma_{14} & \gamma_{15} & \gamma_{16} & \gamma_{17} & \gamma_{18}\\ \gamma_{21} & \gamma_{22} & \gamma_{23} & \gamma_{24} & \gamma_{25} & \gamma_{26} & \gamma_{27} & \gamma_{28}\\ \gamma_{31} & \gamma_{32} & \gamma_{33} & \gamma_{34} & \gamma_{35} & \gamma_{36} & \gamma_{37} & \gamma_{38}\\ \gamma_{41} & \gamma_{42} & \gamma_{43} & \gamma_{44} & \gamma_{45} & \gamma_{46} & \gamma_{47} & \gamma_{48}\\ \gamma_{51} & \gamma_{52} & \gamma_{53} & \gamma_{54} & \gamma_{55} & \gamma_{56} & \gamma_{57} & \gamma_{58}\\ \gamma_{61} & \gamma_{62} & \gamma_{63} & \gamma_{64} & \gamma_{65} & \gamma_{66} & \gamma_{67} & \gamma_{68}\end{bmatrix},$$

$G^*_j(\mu,t)$ is the matrix $\mathfrak{Coef}(\mu)$ when we *replace* its $j^{th}$ column with $[\mathfrak{Const}^*(\mu,t)]^{Tr}$, for $1\le j\le 4$, and $G^*_j(\mu,t)$ is the matrix $\mathfrak{Coef}(\mu)$ when we *replace* its $j^{th}$ column with $[\mathfrak{Const}^*(\mu,t)]^{Tr}$, for $5\le j\le 8$. Here $[\mathfrak{Const}^*(\mu,t)]^{Tr}$ is the transpose of the row

$$\mathfrak{Const}^*(\mu,t):=[J^*_R(\mu,t)+J_R(\mu,t)\quad J^*_L(\mu,t)+J_L(\mu,t)\quad \tilde{I}_1(\omega_0,t)\quad \tilde{I}_2(\omega_0,t)\quad \tilde{I}_3(\omega_0,t)\quad \tilde{I}_4(\omega_0,t)\quad \tilde{I}_5(\omega_0,t)\quad \tilde{I}_6(\omega_0,t)]\,.$$

Substituting (5.5) in (4.3) and (5.6) in (4.6), we find

(5.7) $$u^R(x,t)=\frac{1}{2\pi}\int_{-\infty}^{\infty}\hat{u}_0^R(\lambda)e^{i\lambda x-\omega_R(\lambda)t}d\lambda$$
$$-\frac{1}{2\pi}\int_{\mu\in C(i,\varepsilon)}\frac{i\varphi_R(\mu)e^{i\varphi_R(\mu)x-\omega_0(\mu)t}}{\{1+\alpha_R[\varphi_R(\mu)]^2\}\det[\mathfrak{Coef}(\mu)]}[\alpha_R\det G^*_3(\mu,t)+\beta_R\det G^*_1(\mu,t)]\frac{d[\varphi_R(\mu)]}{d\mu}d\mu$$
$$-\frac{1}{2\pi}\int_{\mu\in C(i,\varepsilon)}\frac{e^{i\varphi_R(\mu)x-\omega_0(\mu)t}}{\{1+\alpha_R[\varphi_R(\mu)]^2\}\det[\mathfrak{Coef}(\mu)]}[\alpha_R\det G^*_4(\mu,t)+\beta_R\det G^*_2(\mu,t)]\frac{d[\varphi_R(\mu)]}{d\mu}d\mu$$
$$+\frac{1}{2\pi}\int_{-\infty}^{\infty}\frac{e^{i\lambda x-\omega_R(\lambda)t}}{1+\alpha_R\lambda^2}\tilde{\hat{f}}_R(\lambda,\omega_R,t)d\lambda\,,\ x>0\,,\ t>0\,,$$

and

(5.8) $$u^L(x,t)=\frac{1}{2\pi}\int_{-\infty}^{\infty}\hat{u}_0^L(\lambda)e^{i\lambda x-\omega_L(\lambda)t}d\lambda$$
$$+\frac{1}{2\pi}\int_{\mu\in C(i,\varepsilon)}\frac{i\varphi_L(\mu)e^{i\varphi_L(\mu)x-\omega_0(\mu)t}}{\{1+\alpha_L[\varphi_L(\mu)]^2\}\det[\mathfrak{Coef}(\mu)]}[\alpha_L\det G_7^*(\mu,t)+\beta_L\det G_5^*(\mu,t)]\frac{d[\varphi_L(\mu)]}{d\mu}d\mu$$
$$+\frac{1}{2\pi}\int_{\mu\in C(i,\varepsilon)}\frac{e^{i\varphi_L(\mu)x-\omega_0(\mu)t}}{\{1+\alpha_L[\varphi_L(\mu)]^2\}\det[\mathfrak{Coef}(\mu)]}[\alpha_L\det G_8^*(\mu,t)+\beta_L\det G_6^*(\mu,t)]\frac{d[\varphi_L(\mu)]}{d\mu}d\mu$$
$$+\frac{1}{2\pi}\int_{-\infty}^{\infty}\frac{e^{i\lambda x-\omega_L(\lambda)t}}{1+\alpha_L\lambda^2}\tilde{\hat{f}}_L(\lambda,\omega_L,t)d\lambda\,,\ x<0\,,\ t>0\,.$$

For the statement of the following proposition, we consider the matrices $G_j(\mu,t)$, which result from the matrices $G_j^*(\mu,t)$, $1\le j\le 8$, by eliminating the terms $J_R^*(\mu,t)$ and $J_L^*(\mu,t)$.

**Proposition** *Suppose that* $\det[\mathfrak{Coef}(i)]\neq 0$ *and*

$$\beta_R\Gamma_{345678}-\alpha_R\Gamma_{235678}=0\,,\ \Gamma_{245678}=0\,,\ \beta_R\Gamma_{345678}-\alpha_R\Gamma_{145678}=0\,,\ \Gamma_{135678}=0\,,$$
$$\beta_R\Gamma_{145678}-\alpha_R\Gamma_{125678}=0\,,\ \beta_R\Gamma_{235678}-\alpha_R\Gamma_{235678}=0\,,$$
$$\beta_L\Gamma_{123478}-\alpha_L\Gamma_{123467}=0\,,\ \Gamma_{123468}=0\,,\ \beta_L\Gamma_{123478}-\alpha_L\Gamma_{123467}=0\,,\ \Gamma_{123457}=0\,,$$
$$\beta_L\Gamma_{123458}-\alpha_L\Gamma_{123456}=0\,,\ \beta_L\Gamma_{123467}-\alpha_L\Gamma_{123456}=0\,.$$

*Then the parts of the integrals in (5.7) and (5.8), taken on* $C(i,\varepsilon)$, *which contain the terms* $J_R^*(\mu,t)$ *and* $J_L^*(\mu,t)$, *vanish, and, therefore, (5.7) and (5.8) hold when we replace* $G_j^*(\mu,t)$ *by* $G_j(\mu,t)$, *for* $1\le j\le 8$.
*Thus, the UTM solution of Problem 3 has the following integral representation:*

(5.9) $$u^R(x,t)=\frac{1}{2\pi}\int_{-\infty}^{\infty}\hat{u}_0^R(\lambda)e^{i\lambda x-\omega_R(\lambda)t}d\lambda$$
$$-\frac{1}{2\pi}\int_{\lambda\in\mathcal{K}(i/\sqrt{\alpha})}\frac{i\lambda e^{i\lambda x-\omega_R(\lambda)t}}{(1+\alpha_R\lambda^2)\det[\mathfrak{Coef}(\psi_R(\lambda))]}[\alpha_R\det G_3(\psi_R(\lambda),t)+\beta_R\det G_1(\psi_R(\lambda),t)]d\lambda$$
$$-\frac{1}{2\pi}\int_{\lambda\in\mathcal{K}(i/\sqrt{\alpha})}\frac{e^{i\lambda x-\omega_R(\lambda)t}}{(1+\alpha_R\lambda^2)\det[\mathfrak{Coef}(\psi_R(\lambda))]}[\alpha_R\det G_4(\psi_R(\lambda),t)+\beta_R\det G_2(\psi_R(\lambda),t)]d\lambda$$
$$+\frac{1}{2\pi}\int_{-\infty}^{\infty}\frac{e^{i\lambda x-\omega_R(\lambda)t}}{1+\alpha_R\lambda^2}\tilde{\hat{f}}_R(\lambda,\omega_R,t)d\lambda\,,\ x>0\,,\ t>0\,,$$

and

(5.10) $$u^L(x,t)=\frac{1}{2\pi}\int_{-\infty}^{\infty}\hat{u}_0^L(\lambda)e^{i\lambda x-\omega_L(\lambda)t}d\lambda$$
$$+\frac{1}{2\pi}\int_{\lambda\in\mathcal{K}(-i/\sqrt{\alpha})}\frac{i\lambda e^{i\lambda x-\omega_L(\lambda)t}}{(1+\alpha_L\lambda^2)\det[\mathfrak{Coef}(\psi_L(\lambda))]}[\alpha_L\det G_7(\psi_L(\lambda),t)+\beta_L\det G_5(\psi_L(\lambda),t)]d\lambda$$
$$+\frac{1}{2\pi}\int_{\lambda\in\mathcal{K}(-i/\sqrt{\alpha})}\frac{e^{i\lambda x-\omega_L(\lambda)t}}{(1+\alpha_L\lambda^2)\det[\mathfrak{Coef}(\psi_L(\lambda))]}[\alpha_L\det G_8(\psi_L(\lambda),t)+\beta_L\det G_6(\psi_L(\lambda),t)]d\lambda$$
$$+\frac{1}{2\pi}\int_{-\infty}^{\infty}\frac{e^{i\lambda x-\omega_L(\lambda)t}}{1+\alpha_L\lambda^2}\tilde{\hat{f}}_L(\lambda,\omega_L,t)d\lambda\,,\ x<0\,,\ t>0\,.$$

*Notation* If $1 \le j_1 < j_2 < j_3 < j_4 < j_5 < j_6 \le 8$, $\Gamma_{j_1 j_2 j_3 j_4 j_5 j_6}$ denotes the determinant of the $6\times 6$ matrix whose columns are the $j_1^{th}, j_2^{th}, j_3^{th}, j_4^{th}, j_5^{th}, j_6^{th}$ columns of the following matrix:

$$\begin{bmatrix} \gamma_{11} & \gamma_{12} & \gamma_{13} & \gamma_{14} & \gamma_{15} & \gamma_{16} & \gamma_{17} & \gamma_{18} \\ \gamma_{21} & \gamma_{22} & \gamma_{23} & \gamma_{24} & \gamma_{25} & \gamma_{26} & \gamma_{27} & \gamma_{28} \\ \gamma_{31} & \gamma_{32} & \gamma_{33} & \gamma_{34} & \gamma_{35} & \gamma_{36} & \gamma_{37} & \gamma_{38} \\ \gamma_{41} & \gamma_{42} & \gamma_{43} & \gamma_{44} & \gamma_{45} & \gamma_{46} & \gamma_{47} & \gamma_{48} \\ \gamma_{51} & \gamma_{52} & \gamma_{53} & \gamma_{54} & \gamma_{55} & \gamma_{56} & \gamma_{57} & \gamma_{58} \\ \gamma_{61} & \gamma_{62} & \gamma_{63} & \gamma_{64} & \gamma_{65} & \gamma_{66} & \gamma_{67} & \gamma_{68} \end{bmatrix}.$$

## 6. The problem 4

***Derivation of the preliminary solution*** The global relation for *Problem 4* [(1.1), (1.14), (1.15)] is given by

$$(6.1)\quad \hat{u}(\lambda,t)=\hat{u}_0(\lambda)e^{-\omega(\lambda)t}-\frac{e^{-\omega(\lambda)t}}{1+\alpha\lambda^2}[i\lambda\widetilde{G}_0(\omega,t)+\widetilde{G}_1(\omega,t)]$$
$$+\frac{e^{-\omega(\lambda)t}}{1+\alpha\lambda^2}e^{-i\lambda l}[i\lambda\widetilde{H}_0(\omega,t)+\widetilde{H}_1(\omega,t)]+\frac{e^{-\omega(\lambda)t}}{1+\alpha\lambda^2}\tilde{\hat{f}}(\lambda,\omega,t),\ \lambda\in\mathbb{C},$$

where

$$\widetilde{G}_0(\omega,t)=\alpha\widetilde{g}'_0(\omega,t)+\beta\widetilde{g}_0(\omega,t),\ \widetilde{G}_1(\omega,t)=\alpha\widetilde{g}'_1(\omega,t)+\beta\widetilde{g}_1(\omega,t),$$
$$\widetilde{H}_0(\omega,t)=\alpha\widetilde{h}'_0(\omega,t)+\beta\widetilde{h}_0(\omega,t),\ \widetilde{H}_1(\omega,t)=\alpha\widetilde{h}'_0(\omega,t)+\beta\widetilde{h}_0(\omega,t),$$
$$g_0(t)=u(0,t),\ g'_0(t)=u_t(0,t),\ g_1(t)=u_x(0,t)\ \text{and}\ g'_1(t)=u_{xt}(0,t),$$
$$g_0(t)=u(l,t),\ h'_0(t):=u_t(l,t),\ h_1(t)=u_x(l,t)\ \text{and}\ h'_1(t):=u_{xt}(l,t).$$

In this section

$$\hat{u}_0(\lambda)=\int_0^l e^{-i\lambda x}u_0(x)dx,\ \hat{u}(\lambda,t)=\int_0^l e^{-i\lambda x}u(x,t)dx,$$
$$\hat{f}(\lambda,t)=\int_0^l e^{-i\lambda x}f(x,t)dx\ \text{ and }\ \tilde{\hat{f}}(\lambda,\omega,t)=\int_0^l e^{\omega(\tau)}\hat{f}(x,\tau)d\tau\ \ (\lambda\in\mathbb{C}).$$

Integrating (6.1), we obtain

$$(6.2)\quad u(x,t)=\frac{1}{2\pi}\int_{-\infty}^{\infty}e^{i\lambda x-\omega(\lambda)t}\hat{u}_0(\lambda)d\lambda-\frac{1}{2\pi}\int_{-\infty}^{\infty}\frac{e^{i\lambda x-\omega(\lambda)t}}{1+\alpha\lambda^2}[i\lambda\widetilde{G}_0(\omega,t)+\widetilde{G}_1(\omega,t)]d\lambda$$
$$+\frac{1}{2\pi}\int_{-\infty}^{\infty}\frac{e^{i\lambda x-\omega(\lambda)t}}{1+\alpha\lambda^2}e^{-i\lambda l}[i\lambda\widetilde{H}_0(\omega,t)+\widetilde{H}_1(\omega,t)]d\lambda+\frac{1}{2\pi}\int_{-\infty}^{\infty}\frac{e^{i\lambda x-\omega(\lambda)t}}{1+\alpha\lambda^2}\tilde{\hat{f}}(\lambda,\omega,t)d\lambda,\ 0<x<l,\ t>0.$$

The integrals in (6.2), involving the terms $i\lambda\widetilde{G}_0(\omega,t)+\widetilde{G}_1(\omega,t)$ and $i\lambda\widetilde{H}_0(\omega,t)+\widetilde{H}_1(\omega,t)$, can be transformed, leading to the *preliminary* representation of the solution:

$$(6.3)\quad u(x,t)=\frac{1}{2\pi}\int_{-\infty}^{\infty}e^{i\lambda x-\omega(\lambda)t}\hat{u}_0(\lambda)d\lambda-\frac{1}{2\pi}\int_{\mathcal{K}(i/\sqrt{\alpha})}\frac{e^{i\lambda x-\omega(\lambda)t}}{1+\alpha\lambda^2}[i\lambda\widetilde{G}_0(\omega,t)+\widetilde{G}_1(\omega,t)]d\lambda$$
$$+\frac{1}{2\pi}\int_{\mathcal{K}^*(-i/\sqrt{\alpha})}\frac{e^{i\lambda x-\omega(\lambda)t}}{1+\alpha\lambda^2}e^{-i\lambda l}[i\lambda\widetilde{H}_0(\omega,t)+\widetilde{H}_1(\omega,t)]d\lambda+\frac{1}{2\pi}\int_{-\infty}^{\infty}\frac{e^{i\lambda x-\omega(\lambda)t}}{1+\alpha\lambda^2}\tilde{\hat{f}}(\lambda,\omega,t)d\lambda,\ 0<x<l,\ t>0,$$

where $\mathcal{K}(i/\sqrt{\alpha})$ is as in fig 1 and $\mathcal{K}^*(-i/\sqrt{\alpha})$ is a sufficiently small simple closed curve surrounding the point $-i/\sqrt{\alpha}$, traversed clockwise.

***Elimination of the unknown quantities*** The global relation (6.1) can be written in the form

(6.4) $$[i\lambda\widetilde{G}_0(\omega,t)+\widetilde{G}_1(\omega,t)]-e^{-i\lambda l}[i\lambda\widetilde{H}_0(\omega,t)+\widetilde{H}_1(\omega,t)]=J(\lambda,t)-J^*(\lambda,t),\ \lambda\in\mathbb{C},$$

where

(6.5) $$J(\lambda,t)=(1+\alpha\lambda^2)\hat{u}_0(\lambda)+\tilde{\hat{f}}(\lambda,\omega,t)\ \text{ and }\ J^*(\lambda,t)=(1+\alpha\lambda^2)e^{\omega(\lambda)t}\hat{u}(\lambda,t).$$

Using the transformation $\lambda\to-\lambda$ in the global relation (4.4), we obtain

(6.6) $$[-i\lambda\widetilde{G}_0(\omega,t)+\widetilde{G}_1(\omega,t)]-e^{-i\lambda l}[-i\lambda\widetilde{H}_0(\omega,t)+\widetilde{H}_1(\omega,t)]=J(-\lambda,t)-J^*(-\lambda,t),\ \lambda\in\mathbb{C}.$$

On the other hand, it follows from the boundary conditions (1.15) that

(6.7) $$\widetilde{G}_1(\omega,t)-c_0\widetilde{g}_0(\omega,t)=\widetilde{\zeta}_0(\omega,t)\ \text{ and }\ \widetilde{H}_1(\omega,t)-c_l\widetilde{h}_0(\omega,t)=\widetilde{\zeta}_l(\omega,t).$$

From the relations

$$\widetilde{g}'_0(\omega,t)=\int_0^t e^{\omega(\lambda)\tau}g_0'(\tau)d\tau=-\omega\widetilde{g}_0(\omega,t)+e^{\omega(\lambda)t}u(0,t)-u_0(0)$$

and

$$\widetilde{h}'_0(\omega,t)=\int_0^t e^{\omega(\lambda)\tau}h_0'(\tau)d\tau=-\omega\widetilde{h}_0(\omega,t)+e^{\omega(\lambda)t}u(l,t)-u_0(l),$$

we obtain

$$\widetilde{G}_0(\omega,t)=(\beta-\alpha\omega)\widetilde{g}_0(\omega,t)+U_0(\omega,t)\ \text{ and }\ \widetilde{H}_0(\omega,t)=(\beta-\alpha\omega)\widetilde{h}_0(\omega,t)+U_l(\omega,t),$$

where

(6.8) $$U_0(\omega,t)=\alpha[e^{\omega t}u(0,t)-u_0(0)]\ \text{ and }\ U_l(\omega,t)=\alpha[e^{\omega t}u(l,t)-u_0(l)].$$

Therefore, the equations (6.4), (6.6) and (6.7) give the following system of equations for $\widetilde{g}_0$ and $\widetilde{h}_0$:

(6.9) $$\begin{bmatrix}\kappa_0(\lambda) & -e^{-i\lambda l}\kappa_l(\lambda)\\ \kappa_0(-\lambda) & -e^{-i\lambda l}\kappa_l(-\lambda)\end{bmatrix}\begin{bmatrix}\widetilde{g}_0(\omega,t)\\ \widetilde{h}_0(\omega,t)\end{bmatrix}$$
$$=\begin{bmatrix}J(\lambda,t)-J^*(\lambda,t)-i\lambda U_0(\omega,t)-\widetilde{\zeta}_0(\omega,t)+e^{-i\lambda l}i\lambda U_l(\omega,t)+e^{-i\lambda l}\widetilde{\zeta}_l(\omega,t)\\ J(-\lambda,t)-J^*(-\lambda,t)+i\lambda U_0(\omega,t)-\widetilde{\zeta}_0(\omega,t)+e^{i\lambda l}i\lambda U_l(\omega,t)+e^{i\lambda l}\widetilde{\zeta}_l(\omega,t)\end{bmatrix},$$

where

$$\kappa_0(\lambda)=i\lambda(\beta-\alpha\omega)+c_0\ \text{ and }\ \kappa_l(\lambda)=i\lambda(\beta-\alpha\omega)+c_l.$$

The determinant of the matrix of the coefficients of the unknowns $\widetilde{g}_0$ and $\widetilde{h}_0$ is

$$\mathcal{D}(\lambda)=\det\begin{bmatrix}\kappa_0(\lambda) & -e^{-i\lambda l}\kappa_l(\lambda)\\ \kappa_0(-\lambda) & -e^{i\lambda l}\kappa_l(-\lambda)\end{bmatrix}=e^{-i\lambda l}\kappa_0(-\lambda)\kappa_l(\lambda)-e^{i\lambda l}\kappa_0(\lambda)\kappa_l(-\lambda).$$

Therefore, working with $\mathcal{D}(\lambda)\neq 0$ and solving the system (6.9), we find that $\widetilde{g}_0$ and $\widetilde{h}_0$ are given by the formulas:

$$\widetilde{g}_0(\omega,t)=\frac{1}{\mathcal{D}(\lambda)}\det\begin{bmatrix}J(\lambda,t)-J^*(\lambda,t)-i\lambda U_0(\omega,t)-\widetilde{\zeta}_0(\omega,t)+e^{-i\lambda l}i\lambda U_l(\omega,t)+e^{-i\lambda l}\widetilde{\zeta}_l(\omega,t) & -e^{-i\lambda l}\kappa_l(\lambda)\\ J(-\lambda,t)-J^*(-\lambda,t)+i\lambda U_0(\omega,t)-\widetilde{\zeta}_0(\omega,t)+e^{i\lambda l}i\lambda U_l(\omega,t)+e^{i\lambda l}\widetilde{\zeta}_l(\omega,t) & -e^{i\lambda l}\kappa_l(-\lambda)\end{bmatrix}$$

and

$$\widetilde{h}_0(\omega,t)=\frac{1}{\mathcal{D}(\lambda)}\det\begin{bmatrix}\kappa_0(\lambda) & J(\lambda,t)-J^*(\lambda,t)-i\lambda U_0(\omega,t)-\widetilde{\zeta}_0(\omega,t)+e^{-i\lambda l}i\lambda U_l(\omega,t)+e^{-i\lambda l}\widetilde{\zeta}_l(\omega,t)\\ \kappa_0(-\lambda) & J(-\lambda,t)-J^*(-\lambda,t)+i\lambda U_0(\omega,t)-\widetilde{\zeta}_0(\omega,t)+e^{i\lambda l}i\lambda U_l(\omega,t)+e^{i\lambda l}\widetilde{\zeta}_l(\omega,t)\end{bmatrix}.$$

These, in turn, give

(6.10)
$$\begin{aligned}i\lambda\widetilde{G}_0(\omega,t)+\widetilde{G}_1(\omega,t)&=\kappa_0(\lambda)\widetilde{g}_0(\omega,t)+i\lambda U_0(\omega,t)+\widetilde{\zeta}_0(\omega,t)\\ &=[\mathcal{D}(\lambda)]^{-1}\kappa_0(\lambda)[-e^{i\lambda l}\kappa_l(-\lambda)J(\lambda,t)+e^{-i\lambda l}\kappa_l(\lambda)J(-\lambda,t)\\ &\quad+[\mathcal{D}(\lambda)]^{-1}\kappa_0(\lambda)[e^{i\lambda l}\kappa_l(-\lambda)J^*(\lambda,t)-e^{-i\lambda l}\kappa_l(\lambda)J^*(-\lambda,t)]\\ &\quad+[\mathcal{D}(\lambda)]^{-1}e^{-i\lambda l}2i\lambda\kappa_l(\lambda)[c_0U_0(\omega,t)-(\beta-\alpha\omega)\widetilde{\zeta}_0(\omega,t)]\\ &\quad+[\mathcal{D}(\lambda)]^{-1}2i\lambda\kappa_0(\lambda)[c_lU_l(\omega,t)+(\beta-\alpha\omega)\widetilde{\zeta}_l(\omega,t)],\end{aligned}$$

and

(6.11)
$$\begin{aligned}i\lambda\widetilde{H}_0(\omega,t)+\widetilde{H}_1(\omega,t)&=e^{-i\lambda l}[\kappa_l(\lambda)\widetilde{h}_0(\omega,t)+i\lambda U_l(\omega,t)+\widetilde{\zeta}_l(\omega,t)]\\ &=[\mathcal{D}(\lambda)]^{-1}\kappa_l(\lambda)e^{-i\lambda l}[-\kappa_0(-\lambda)J(\lambda,t)+\kappa_0(\lambda)J(-\lambda,t)]\\ &\quad+[\mathcal{D}(\lambda)]^{-1}\kappa_l(\lambda)e^{-i\lambda l}[\kappa_0(-\lambda)J^*(\lambda,t)-\kappa_0(\lambda)J^*(-\lambda,t)]\\ &\quad+[\mathcal{D}(\lambda)]^{-1}e^{-i\lambda l}2i\lambda\kappa_l(\lambda)[c_0U_0(\omega,t)-(\beta-\alpha\omega)\widetilde{\zeta}_0(\omega,t)]\\ &\quad+[\mathcal{D}(\lambda)]^{-1}2i\lambda\kappa_0(\lambda)[c_lU_l(\omega,t)+(\beta-\alpha\omega)\widetilde{\zeta}_l(\omega,t)].\end{aligned}$$

From the asymptotic behavior

$$\omega(\lambda)\sim-\frac{\beta}{\alpha}\frac{1}{1+\alpha\lambda^2}\quad\text{as}\quad\lambda\to\pm\frac{i}{\sqrt{\alpha}},$$

it follows that

$$\kappa_0(\lambda),\ \kappa_l(\lambda)\sim i\lambda\beta\frac{1}{1+\alpha\lambda^2}\quad\text{as}\quad\lambda\to\pm\frac{i}{\sqrt{\alpha}}$$

and

$$\frac{1}{\mathcal{D}(\lambda)}\sim\frac{\alpha}{\beta^2}\frac{(1+\alpha\lambda^2)^2}{e^{i\lambda l}-e^{-i\lambda l}}\equiv\frac{\alpha}{\beta^2}\frac{(1+\alpha\lambda^2)^2}{2\sin(\lambda l)}\quad\text{as}\quad\lambda\to\pm\frac{i}{\sqrt{\alpha}}.$$

Therefore, using the expression for $J^*$ (see equation (4.5)), we obtain

(6.12)
$$\begin{aligned}&\int_{\mathcal{K}(i/\sqrt{\alpha})}\frac{e^{i\lambda x-\omega(\lambda)t}}{1+\alpha\lambda^2}\frac{\kappa_0(\lambda)[e^{i\lambda l}\kappa_l(-\lambda)J^*(\lambda,t)-e^{-i\lambda l}\kappa_l(\lambda)J^*(-\lambda,t)]}{\mathcal{D}(\lambda)}d\lambda\\ &\qquad=\int_{\mathcal{K}(i/\sqrt{\alpha})}e^{i\lambda x}\frac{\kappa_0(\lambda)[e^{i\lambda l}\kappa_l(-\lambda)\hat{u}(\lambda,t)-e^{-i\lambda l}\kappa_l(\lambda)\hat{u}(-\lambda,t)]}{\mathcal{D}(\lambda)}d\lambda=0\end{aligned}$$

and

(6.13)
$$\begin{aligned}&\int_{\mathcal{K}^*(-i/\sqrt{\alpha})}\frac{e^{i\lambda x-\omega(\lambda)t}}{1+\alpha\lambda^2}\frac{\kappa_l(\lambda)e^{-i\lambda l}[\kappa_0(-\lambda)J^*(\lambda,t)-\kappa_0(\lambda)J^*(-\lambda,t)]}{\mathcal{D}(\lambda)}d\lambda\\ &\qquad\int_{\mathcal{K}^*(-i/\sqrt{\alpha})}e^{i\lambda x}\frac{\kappa_l(\lambda)e^{-i\lambda l}[\kappa_0(-\lambda)\hat{u}(\lambda,t)-\kappa_0(\lambda)\hat{u}(-\lambda,t)]}{\mathcal{D}(\lambda)}d\lambda=0,\end{aligned}$$

because the integrands are bounded in the vicinities of the points $\pm i/\sqrt{\alpha}$.

In the same way, using the equations (6.8), we obtain

$$(6.14)\quad \int_{\mathcal{K}(i/\sqrt{\alpha})} \frac{e^{i\lambda x-\omega(\lambda)t}}{1+\alpha\lambda^2} 2i\lambda \frac{e^{-i\lambda l}\kappa_l(\lambda)c_0U_0(\omega,t)+\kappa_0(\lambda)c_lU_l(\omega,t)}{\mathcal{D}(\lambda)}d\lambda$$
$$=\alpha \int_{\mathcal{K}(i/\sqrt{\alpha})} \frac{e^{i\lambda x}}{1+\alpha\lambda^2} 2i\lambda \frac{e^{-i\lambda l}\kappa_l(\lambda)c_0u(0,t)+\kappa_0(\lambda)c_lu(l,t)}{\mathcal{D}(\lambda)}d\lambda$$
$$-\alpha \int_{\mathcal{K}(i/\sqrt{\alpha})} \frac{e^{i\lambda x-\omega(\lambda)t}}{1+\alpha\lambda^2} 2i\lambda \frac{e^{-i\lambda l}\kappa_l(\lambda)c_0u_0(0)+\kappa_0(\lambda)c_lu_0(l)}{\mathcal{D}(\lambda)}d\lambda$$
$$=-\alpha \int_{\mathcal{K}(i/\sqrt{\alpha})} \frac{e^{i\lambda x-\omega(\lambda)t}}{1+\alpha\lambda^2} 2i\lambda \frac{e^{-i\lambda l}\kappa_l(\lambda)c_0u_0(0)+\kappa_0(\lambda)c_lu_0(l)}{\mathcal{D}(\lambda)}d\lambda$$

and, similarly,

$$(6.15)\quad \int_{\mathcal{K}^*(-i/\sqrt{\alpha})} \frac{e^{i\lambda x-\omega(\lambda)t}}{1+\alpha\lambda^2} 2i\lambda \frac{e^{-i\lambda l}\kappa_l(\lambda)c_0U_0(\omega,t)+\kappa_0(\lambda)c_lU_l(\omega,t)}{\mathcal{D}(\lambda)}d\lambda$$
$$=-\alpha \int_{\mathcal{K}^*(-i/\sqrt{\alpha})} \frac{e^{i\lambda x-\omega(\lambda)t}}{1+\alpha\lambda^2} 2i\lambda \frac{e^{-i\lambda l}\kappa_l(\lambda)c_0u_0(0)+\kappa_0(\lambda)c_lu_0(l)}{\mathcal{D}(\lambda)}d\lambda .$$

Finally, in order to obtain the *genuine* representation of the solution, we use the preliminary representation (6.3), the formulae (6.10) and (6.11), the equation (6.5) for $J(\omega,t)$ and the equalities (6.12), (6.13), (6.14) and (6.15). As a result, the solution to *Problem 4* is given by the formula:

$$(6.16)\quad u(x,t)=\frac{1}{2\pi}\int_{-\infty}^{\infty} e^{i\lambda x-\omega(\lambda)t}\hat{u}_0(\lambda)d\lambda+\frac{1}{2\pi}\int_{-\infty}^{\infty}\frac{e^{i\lambda x-\omega(\lambda)t}}{1+\alpha\lambda^2}\tilde{\hat{f}}(\lambda,\omega,t)d\lambda$$
$$-\frac{1}{2\pi}\int_{\mathcal{K}(i/\sqrt{\alpha})} e^{i\lambda x-\omega(\lambda)t}\frac{\kappa_0(\lambda)[-e^{i\lambda l}\kappa_l(-\lambda)\hat{u}_0(\lambda)+e^{-i\lambda l}\kappa_l(\lambda)\hat{u}_0(-\lambda)]}{\mathcal{D}(\lambda)}d\lambda$$
$$-\frac{1}{2\pi}\int_{\mathcal{K}(i/\sqrt{\alpha})} \frac{e^{i\lambda x-\omega(\lambda)t}}{1+\alpha\lambda^2}\frac{\kappa_0(\lambda)[-e^{i\lambda l}\kappa_l(-\lambda)\tilde{\hat{f}}(\lambda,\omega,t)+e^{-i\lambda l}\kappa_l(\lambda)\tilde{\hat{f}}(-\lambda,\omega,t)]}{\mathcal{D}(\lambda)}d\lambda$$
$$-\frac{i}{\pi}\int_{\mathcal{K}(i/\sqrt{\alpha})} \frac{e^{i\lambda x-\omega(\lambda)t}}{1+\alpha\lambda^2}\frac{\lambda e^{-i\lambda l}\kappa_l(\lambda)[-\alpha c_0u_0(0)-(\beta-\alpha\omega)\tilde{\zeta}_0(\omega,t)]}{\mathcal{D}(\lambda)}d\lambda$$
$$-\frac{i}{\pi}\int_{\mathcal{K}(i/\sqrt{\alpha})} \frac{e^{i\lambda x-\omega(\lambda)t}}{1+\alpha\lambda^2}\frac{\lambda\kappa_0(\lambda)[-\alpha c_lu_0(l)-(\beta-\alpha\omega)\tilde{\zeta}_l(\omega,t)]}{\mathcal{D}(\lambda)}d\lambda$$
$$-\frac{1}{2\pi}\int_{\mathcal{K}^*(-i/\sqrt{\alpha})} e^{i\lambda x-\omega(\lambda)t}\frac{\kappa_l(\lambda)e^{-i\lambda l}[-\kappa_0(-\lambda)\hat{u}_0(\lambda)+\kappa_0(\lambda)\hat{u}_0(-\lambda)}{\mathcal{D}(\lambda)}d\lambda$$
$$-\frac{1}{2\pi}\int_{\mathcal{K}^*(-i/\sqrt{\alpha})} \frac{e^{i\lambda x-\omega(\lambda)t}}{1+\alpha\lambda^2}\frac{\kappa_l(\lambda)e^{-i\lambda l}[-\kappa_0(-\lambda)\tilde{\hat{f}}(\lambda,\omega,t)+\kappa_0(\lambda)\tilde{\hat{f}}(-\lambda,\omega,t)}{\mathcal{D}(\lambda)}d\lambda$$
$$-\frac{i}{\pi}\int_{\mathcal{K}^*(-i/\sqrt{\alpha})} \frac{e^{i\lambda x-\omega(\lambda)t}}{1+\alpha\lambda^2}\frac{\lambda e^{-i\lambda l}\kappa_l(\lambda)[-\alpha c_0u_0(0)-(\beta-\alpha\omega)\tilde{\zeta}_0(\omega,t)]}{\mathcal{D}(\lambda)}d\lambda$$
$$-\frac{i}{\pi}\int_{\mathcal{K}^*(-i/\sqrt{\alpha})} \frac{e^{i\lambda x-\omega(\lambda)t}}{1+\alpha\lambda^2}\frac{\lambda\kappa_0(\lambda)[-\alpha c_lu_0(l)-(\beta-\alpha\omega)\tilde{\zeta}_l(\omega,t)]}{\mathcal{D}(\lambda)}d\lambda ,\ \text{for } 0<x<l,\ t>0.$$

We have thus obtained novel solution formulae for interface as well as initial-boundary-value problems for the Barenblatt pseudo-parabolic equation, via suitable extension and implementation of the UTM methodology.

## Declarations

*Ethical Approval*: Not applicable. *Competing Interests*: Not applicable. *Availability of data and materials*: Not applicable.

## Acknowledgment

A.C. and E.A. are grateful to the Stavros Niarchos Foundation (SNF) and the Hellenic Foundation for Research and Innovation (H.F.R.I.) which co-funded this research under the 5th Call of "Science and Society" Action, "*Always strive for excellence* — Theodoros Papazoglou" (Project Number: 28085). A.C. would like to express his thankfulness to Prof. Avraam A. Konstantinidis (Laboratory of Mechanics and Materials, Aristotle University of Thessaloniki) for academic support, to Prof. Harris Wong (Department of Mechanical and Industrial Engineering, Louisiana State University) for useful comments and encouraging discussions, as well as to Prof. Martin Z. Bazant (Department of Mathematics and Department of Chemical Engineering, Massachusetts Institute of Technology) and Prof. Youcef Mammeri (Institut Camille Jordan CNRS, Université Jean Monnet) for valuable private communications pertaining to PDE-based modelling in battery research.